\documentclass[11pt]{amsart} 
\usepackage{verbatim, latexsym, amssymb, amsmath,color}
\usepackage{epsfig}

\def\N{\mathbb N}
\def\Z{\mathbb Z}

\def\d{\mathrm{div}}

\theoremstyle{remark}

\theoremstyle{definition}

\title{Morse index and multiplicity of min-max  minimal hypersurfaces}
\author{Fernando C. Marques and Andr\'e Neves}
\address{Princeton University \\ Fine Hall \\ Princeton NJ 08544 \\ USA}
\email{coda@math.princeton.edu}
\address{Imperial College London\\ Huxley Building \\ 180 Queen's Gate \\ London SW7 2RH \\ United Kingdom}
\email{a.neves@imperial.ac.uk}
\thanks{The first author was partly supported by NSF-DMS-1509027. The second author was partly supported by  ERC-2011-StG-278940 and EPSRC Programme Grant EP/K00865X/1.}

\begin{document}
\maketitle

\begin{abstract}
The Min-max Theory for the area functional, started by Almgren in the  early 1960s and greatly improved by Pitts in 1981, was left incomplete because it gave no  Morse index estimate for the min-max minimal hypersurface.

We advance the theory further and prove the first general Morse index bounds for minimal hypersurfaces produced by it.
We also settle the multiplicity problem for the classical case of one-parameter sweepouts.

 \end{abstract}

\section{Introduction}

The Min-max Theory for the area functional was started by Almgren in the early 1960s \cite{almgren, almgren-varifolds} and greatly improved  by Pitts  in 1981, resulting in  the proof of the  Almgren-Pitts Min-max Theorem \cite{pitts}. The subject can be regarded as a deep higher-dimensional generalization of the study of closed geodesics and 
uses the regularity results of  Schoen-Simon \cite{schoen-simon}.


Despite the many applications of this theory found recently, no general information about the Morse index was known. The difficulty is mainly caused by the weak notions of convergence (in the sense of Geometric Measure Theory) and the possibility of multiplicity. There is no Hilbert space structure or a Palais-Smale condition to check.

  In this paper we prove the first general Morse index bounds for minimal hypersurfaces produced by this theory. Inspired by finite-dimensional Morse theory, one expects that the Morse index of the critical point being produced should be less than or equal to the number of parameters (dimension of the families used in the min-max process). Here we confirm this heuristic for the index of the support of the minimal hypersurface, i.e. the sum of the indices of the components.  

We also prove the first general multiplicity one theorem in min-max theory, in the setting of one-parameter sweepouts. These sweepouts have been  considered since the work of  Birkhoff \cite{birkhoff}.
 We show that for manifolds with generic metrics, the min-max minimal hypersurface either has Morse index one with a multiplicity one  two-sided unstable  component,  or is stable (i.e. has Morse index zero) and contains 
an one-sided component with even multiplicity and unstable  double cover.  We do not know whether the second case can be perturbed away in general, except of course if there are no one-sided closed hypersurfaces like in the case of manifolds $M^{n+1}$ with $H_{n}(M,\mathbb{Z}_2)=0$. Nonetheless, the fact that a two-sided min-max minimal hypersurface must be unstable generically was not known prior to our work.

Motivated by our finding, we propose the  following conjecture for the $k$-parameter min-max setting:
\medskip

{\bf Multiplicity One Conjecture:} {\it For generic metrics on $M^{n+1}$,  $3\leq (n+1)\leq 7$,  two-sided unstable components of  closed minimal hypersurfaces obtained by min-max methods must have multiplicity one.}
\medskip

Theorem \ref{theorem2} of this paper proves this conjecture in the case of min-max with one parameter. 

The one-parameter case in manifolds with positive Ricci curvature had been studied before by Zhou  \cite{zhou, zhou2}, who did important work extending to high dimensions previous work of  the authors \cite{marques-neves-duke}. These results have been improved recently by Ketover and the authors in \cite{ketover-marques-neves}, where new index characterizations and multiplicity one theorems are proven as a consequence of the catenoid estimate. Related results have been proven for the least-area closed minimal hypersurface by Mazet and Rosenberg \cite{mazet-rosenberg} and later by Song \cite{song}.

 \subsection{Remark:} For surfaces the conjecture is not true as it is shown in Aiex \cite{nicolau}. This case is special  because when $n=1$ the Almgren-Pitts Min-max Theorem provides the existence of a stationary geodesic network, not necessarily smooth. \medskip

We now state the main theorems and some applications.

Let $(M^{n+1},g)$ be an $(n+1)$-dimensional closed Riemannian manifold,  with $3\leq (n+1)\leq 7$. 
Let $X$ be a simplicial complex of dimension $k$ and $\Phi:X \rightarrow \mathcal Z_n(M^{n+1};{\bf F};G)$ be a continuous map. Here $G$ is either  the group $\mathbb{Z}$ or $\mathbb{Z}_2$, while $\mathcal Z_n(M;G)$ denotes the space of $n$-dimensional flat chains $T$
 in $M$ with coefficients in $G$ and $\partial T=0$. These are called {\it flat cycles}. The notation $\mathcal Z_n(M;{\bf F};G)$ indicates the space $\mathcal Z_n(M;G)$ endowed with the $F$-metric, to be defined later. Basically this means $\Phi$ is continuous in both the flat and the varifold topologies.
 
We let $\Pi$ be the class of all continuous maps $\Phi':X \rightarrow \mathcal Z_n(M^{n+1};{\bf F};G)$ such that $\Phi$ and $\Phi'$ are homotopic to each other in the flat topology. 
 The {\it width} of $\Pi$ is defined to be the min-max invariant:
 $$
 {\bf L}(\Pi) = \inf_{\Phi' \in \Pi}\sup_{x\in X}\{{\bf M}(\Phi'(x))\},
 $$
 where ${\bf M}(T)$ denotes the mass of $T$ (or $n$-dimensional area).
 
 Given a sequence  $\{\Phi_i\}_{i\in\N}$  of continuous maps from $X$ into $\mathcal Z_n(M;{\bf F};G)$, we set 
$${\bf L}(\{\Phi_i\}_{i\in\N}):=\limsup_{i \rightarrow \infty} \sup_{x\in X} {\bf M}(\Phi_i(x)).$$
When ${\bf L}(\{\Phi_i\}_{i\in\N})={\bf L}(\Pi)$, we say $\{\Phi_i\}_{i\in\N}$ is a {\it min-max sequence} in $\Pi$. 

The Almgren-Pitts min-max theory is actually formulated in terms of sequences of maps defined on the vertices of finer and finer subdivisions of $X$ and that become finer and finer in the mass topology. This imposes quite a bit of notation and technicalities. In Section \ref{minmax.continuous}, we explain how to use the original discretized setting to obtain a min-max theory for maps that are continuous in the ${\bf F}$-topology as above.   A similar construction for the mass topology was done by the authors in \cite{marques-neves-cdm}.


The next theorem establishes general upper Morse index bounds for the min-max minimal hypersurface. It is a consequence of Deformation Theorem A of Section \ref{deformations} and its proof is left to Section \ref{index.bounds}.

\subsection{Theorem}\label{theorem1}{\em Let $(M^{n+1},g)$ be an $(n+1)$-dimensional closed Riemannian manifold,  with $3\leq (n+1)\leq 7$.  There exists an integral stationary varifold $\Sigma \in \mathcal V_n(M)$ with the following properties:
\begin{itemize}
\item $||\Sigma||(M)={\bf L}(\Pi),$
\item the support of $\Sigma$ is a smooth closed embedded hypersurface in $M$,
\item $index(support \, of \, \Sigma) \leq k$.
\end{itemize}
}

\subsection{Simon-Smith variant} For three-dimensional manifolds $M^3$, one can restrict to sweepouts by smooth surfaces (or by level sets of Morse functions) and consider the class $\Pi$ of all sweepouts that differ from a fixed one by the action of ambient isotopies. The isotopies vary smoothly in the parameter space. This is the Simon-Smith variant of Almgren-Pitts theory (see Smith \cite{smith} and a survey by Colding-De Lellis \cite{colding-delellis}). This variant is specially convenient if one wants to control the topology of the min-max minimal surface. Genus bounds were proven recently by Ketover
\cite{ketover} (see also De Lellis-Pellandini \cite{delellis-genus}), as conjectured by Pitts and Rubinstein \cite{pitts-rubinstein}. We note that another variant of min-max theory, based on phase transition partial differential equations, has been recently proposed by Guaraco \cite{guaraco}.

If one follows the proof of Theorem \ref{theorem1}, which is based on Deformation Theorem A, one will see that all deformations are by isotopies. By mollifying these deformations so that they depend smoothly of the parameters, we get that Theorem \ref{theorem1} is also true in the Simon-Smith setting. We leave the details to the reader. As an application, the main result of Ketover-Zhou \cite{ketover-zhou} about the entropy of closed surfaces can be extended to surfaces with genus different than one.

\subsection{Willmore and Freedman-He-Wang conjectures}
 We proved the Willmore conjecture in \cite{marques-neves} by introducing certain five-parameter sweepouts of closed  surfaces in the three-sphere and using min-max theory. Later we proved jointly with I. Agol \cite{agol-marques-neves} the Freedman-He-Wang conjecture about links following a similar strategy. The basic intuition was that the min-max minimal surface produced would have index less than or equal to five (the number of parameters) and then it could be only an equatorial sphere or the Clifford torus by a theorem of Urbano \cite{urbano}. The possibility of an equatorial sphere is ruled out by topological arguments. The index bound was not available at the time of those papers, and we had to proceed indirectly. With Theorem \ref{theorem1} we can argue directly and conclude the min-max surface is a Clifford torus. This simplifies a bit the work involved in the proofs of the main statements of \cite{marques-neves} and \cite{agol-marques-neves}.
 
 \subsection{The space of cycles} The calculation by Almgren \cite{almgren} of the homotopy groups of the space of $n$-cycles mod 2 in $M^{n+1}$ gives that this space is weakly homotopically equivalent to $\mathbb{RP}^\infty$ (see \cite{marques-neves-infinitely}). 
 
 We used this structure in \cite{marques-neves-infinitely} to prove Yau's conjecture \cite{yau1} about the existence of  infinitely many closed minimal hypersurfaces, in the case the manifold has positive Ricci curvature. The minimal hypersurfaces are constructed by applying min-max theory to the class of $k$-sweepouts, which are multiparameter families defined by a cohomological condition, for each $k$. We now explain how to deduce from Theorem \ref{theorem1} that the support of the corresponding min-max minimal hypersurface must have index less than or equal to $k$. There has been considerable activity after our paper \cite{marques-neves-infinitely} regarding compactness properties of minimal hypersurfaces with bounded index (\cite{buzano-sharp}, \cite{carlotto}, \cite{chodosh-ketover-maximo}, \cite{li-zhou}).
 
We have $H^1(\mathcal Z_n(M;\Z_2);\Z_2)=\Z_2=\{0,\bar\lambda\}$ and, for every $k\in\N$, we consider  the set $\mathcal{P}_k$ of all continuous maps $\Phi:X \rightarrow \mathcal Z_n(M;{\bf F};\Z_2)$  such that $\Phi^*(\bar \lambda)^k$ does not vanish in $H^k(X;\Z_2)$, where $X$ is some finite dimensional compact simplicial complex. We defined in \cite{marques-neves-infinitely}
$$\omega_k:=\inf_{\Phi\in\mathcal{P}_k}\sup_{x\in dmn(\Phi)}{\bf M}(\Phi(x)).$$
The first remark is that if $\tilde{\mathcal{P}}_k$ denotes those elements  $\Phi\in \mathcal{P}_k$ whose domain of definition  is  a compact simplicial complex of dimension $k$ then
$$\omega_k=\inf_{\Phi\in\tilde{\mathcal{P}}_k}\sup_{x\in dmn(\Phi)}{\bf M}(\Phi(x)).$$
Indeed, given $\Phi\in \mathcal{P}_k$ consider $X^{(k)}$ the $k$-skeleton of $X:=dmn(\Phi)$. Then $H^k(X,X^{(k)};\Z_2)=0$ and so the long exact cohomology sequence gives that the natural pullback map from $H^k(X;\Z_2)$ into $H^k(X^{(k)};\Z_2)$ is injective, which means that $\Phi_{|X^{(k)}}\in \mathcal{P}_k$ and so $\Phi_{|X^{(k)}}\in \tilde{\mathcal{P}}_k.$

Consider a min-max sequence $\{\Phi_i\}_{i\in\N}\subset   \tilde{\mathcal{P}}_k$ for $\omega_k$. Each $\Phi_i$ induces a homotopy class $\Pi_i$, $i\in\N$. From Theorem \ref{theorem1} we obtain that ${\bf L}(\Pi_i)=||\Sigma_i||(M)$, where the support of $\Sigma_i$ has index $\leq k$. Note that ${\bf L}(\Pi_i)$ tends to $\omega_k$. The Compactness Theorem proven by Sharp \cite{sharp} implies that, after passing to a subsequence, $\{\Sigma_i\}$ converges in the varifold sense to an integral stationary varifold $\Sigma$ with support a minimal embedded hypersurface with  index $\leq k$ and $||\Sigma||(M)=\omega_k$.


\subsection{One-parameter case} Next we assume that $X$ is the interval $[0,1]$, for simplicity. We further assume that the continuous maps satisfy $\Phi(0)=\Phi(1)=0$ and the homotopy class is relative to $\partial [0,1]=\{0,1\}$. The case $X=S^1$ is also considered in the paper.

The metric $g$ is  called {\it bumpy} if there is no smooth immersed minimal hypersurface with a non-trivial Jacobi  field. White showed in \cite{white2, white3} that bumpy metrics are generic in the Baire sense.

The next theorem establishes the instability of the min-max minimal hypersurface, and proves that the multiplicity of an unstable two-sided component must be exactly equal to one. The proof uses the Deformation Theorems B and C and is left to Section \ref{index.bounds}.

\subsection{Theorem}\label{theorem2}{\em  Suppose $(M,g)$ is a bumpy metric, and $X=[0,1]$. For any min-max sequence $\{\Phi_i\}_{i\in\N}$  in $\Pi$, there exists a subsequence 
$i_j \rightarrow \infty$ such that, for some $x_{i_j}\in X$, $\Phi_{i_j}(x_{i_j})$ converges in varifold sense to an integral stationary varifold $\Sigma \in \mathcal V_n(M)$ with $||\Sigma||(M)={\bf L}(\Pi)$ and satisfying:
 \begin{itemize}
 \item[(a)] the support of $\Sigma$ is a smooth closed embedded minimal hypersurface in $M$;
 \end{itemize}
 and either
 \begin{itemize}
 \item[(b1)] $\text{ index\,(support of $\Sigma$)}=1,$ and the unstable component of $\Sigma$ has multiplicity one;
 \end{itemize}
 or
 \begin{itemize}
 \item[(b2)] $\text{ index\,(support of $\Sigma$)}=0,$ and at least one component of $\Sigma$ is one-sided, has even multiplicity and  unstable double cover.
 \end{itemize}
}

\medskip
The theorem has the following  corollary:

\subsection{Corollary}\label{corollary1}{\em  Suppose $(M,g)$ contains no one-sided embedded hypersurfaces, and $X=[0,1]$.
There exists an integral stationary varifold $\Sigma \in \mathcal V_n(M)$ with $||\Sigma||(M)={\bf L}(\Pi)$ and satisfying  the following properties:
\begin{itemize}
\item[(a)] the support of $\Sigma$ is a smooth closed embedded minimal hypersurface in $M$;
\item[(b)]  we have \begin{eqnarray*}
index(support \, of \, \Sigma)  &\leq& 1\\
&\leq& index(support \, of \, \Sigma) +nullity(support \, of \, \Sigma);
\end{eqnarray*}
\item[(c)] any unstable component (necessarily unique)  must have  multiplicity one;
\item[(d)] if the metric $g$ is bumpy and $\{\Phi_i\}_{i\in\N}$ is a min-max sequence, $\Sigma$ can be chosen so that for some subsequence 
$i_j \rightarrow \infty$  there is some $x_{i_j}\in X$ such that $\Phi_{i_j}(x_{i_j})$ converges in varifold sense to $\Sigma$.
\end{itemize}
}
Notice that the topological assumption being required in Corollary \ref{corollary1} is automatically satisfied if $H_{n-1}(M, \mathbb{Z}_2)=0$.

This paper is organized as follows. In Section 2, we set up some basic notation for the rest of the article. In Section 3, we describe how to formulate
a min-max theory that  makes sense in a continuous setting. In Section 4, we define instability of stationary varifolds and derive some consequences.
In Section 5, we prove Deformation Theorems A, B and C which are the fundamental tools in our analysis. In Section 6, we prove the main theorems on the Morse index and multiplicity of min-max minimal hypersurfaces.

\section{Basic notation}\label{gmt}

Let $(M^{n+1},g)$ be an $(n+1)$-dimensional closed Riemannian manifold.  We assume, for convenience, that $(M,g)$ is isometrically embedded in some Euclidean space $\mathbb{R}^L$.

Let $G$ be either the group $\mathbb{Z}$ or $\mathbb{Z}_2$. The spaces we will work with in this paper are:
\begin{itemize}
\item the space ${\bf I}_{n}(M;G)$  of $n$-dimensional  flat chains    in $\mathbb{R}^L$ with coefficients in $G$ and support contained  in $M$; 
\item the space ${\mathcal Z}_n(M;G)$  of flat chains  $T \in {\bf I}_n(M;G)$ with  $\partial T=0$. These are called {\it flat cycles};
\item the space ${\mathcal Z}_n(U;G)$ of flat cycles with support contained in the open set U;
\item the closure $\mathcal{V}_n(M)$, in the weak topology, of the space of $n$-dimensional rectifiable varifolds in $\mathbb{R}^L$ with support contained in $M$. 
\end{itemize}

Given $T\in {\bf I}_n(M;G)$,  we denote by $|T|$ and $||T||$ the integral varifold   and the Radon measure in $M$ associated with $|T|$, respectively;  given $V\in \mathcal{V}_n(M)$, $||V||$ denotes the Radon measure in $M$ associated with $V$.  

 The above spaces come with several relevant metrics. The  {\it flat metric} and the {\it mass} of $T \in {\bf I}_n(M;G)$, denoted by $\mathcal{F}(T)$ and ${\bf M}(T)$,  are defined in \cite{federer}, 
  respectively.  The  ${\bf F}$-{\it metric}  is defined in  {Pitts book} \cite[page 66]{pitts} and   induces the varifold weak topology on $\mathcal{V}_n(M)\cap \{V: ||V||(M) \leq a\}$ for any $a$. We denote by ${\overline{\bf B}^{\bf F}_{\delta}(V)}$ and ${\bf B}^{\bf F}_{\delta}(V)$ the closed and open metric balls, respectively, with radius $\delta$ and center $V \in \mathcal{V}_n(M)$. Similarly, we denote by ${\overline{\bf B}^{\mathcal F}_{\delta}(T)}$ and ${\bf B}^{\mathcal F}_{\delta}(T)$ the corresponding balls with  center $T \in \mathcal{Z}_n(M;G)$ in the flat metric. 
Finally,  the ${\bf F}$-{\it metric} on ${\bf I}_n(M;G)$ is defined by
$$ {\bf F}(S,T)=\mathcal{F}(S-T)+{\bf F}(|S|,|T|).$$

We assume that  ${\bf I}_n(M;G)$ and  ${\mathcal Z}_n(M;G)$ have the topology induced by the flat metric. When endowed with
the topology of the ${\bf F}$-metric or the mass norm, these spaces will be denoted by  ${\bf I}_n(M;{\bf F};G)$, ${\mathcal Z}_n(M;{\bf F};G)$, ${\bf I}_n(M;{\bf M};G)$, ${\mathcal Z}_n(M;{\bf M};G)$, respectively. The space $\mathcal{V}_n(M)$ is considered with the weak topology of varifolds.
Given $\mathcal{A,B}\subset \mathcal{V}_n(M)$, we also define
  $${\bf F}(\mathcal{A},\mathcal{B})=\inf\{{\bf F}(V,W):V\in \mathcal{A}, W\in \mathcal{B}\}.$$

We take $X$ to be a  cubical complex of dimension $k$, i.e. a subcomplex of dimension  $k$ of the  $m$-dimensional cube $I^m=[0,1]^m$ for some $m$. Each $p$-cell of $I^m$ is of the form  $\alpha_1 \otimes \cdots \otimes \alpha_m$, where $\alpha_i \in \{0,1,[0,1]\}$ for every $i$ and $\sum \textrm{dim}(\alpha_i)=p$. By Construction 4.8 of \cite{panov}, any $k$-dimensional simplicial complex is homeomorphic to a  cubical complex of dimension $k$.  If $Y$ is a subcomplex of $X$, we denote by $Y_0$ the set of vertices of $Y$.
  
Throughout this paper we will make use of results contained in \cite{marques-neves} and \cite{marques-neves-infinitely}.  In Section 2 of \cite{marques-neves-infinitely} we described  modifications to the Almgren-Pitts Min-Max Theory so that the $k$-dimensional cube $I^k$ is replaced by $X$ as the parameter space. We refer the reader to Section 2 of \cite{marques-neves-infinitely} for more
details of the Almgren-Pitts  discrete formulation of min-max theory. In that formulation, continuous maps $\Phi$ defined on $X$ are replaced by sequences of maps $\phi_i$ defined on the vertices of finer and finer subdivisions of $X$, so that $M(\phi_i(x)-\phi_i(y))$ becomes arbitrarily small whenever $x,y$ are adjacent vertices  as $i \rightarrow \infty$. 

\section{Min-max theory in continuous setting}\label{minmax.continuous}

As mentioned in the last section, the Almgren-Pitts min-max theory is formulated in terms of sequences of maps defined on the vertices of finer and finer grids. This requires quite a bit of notation and  a lot of technical work. Besides, one works with the mass topology (\cite{pitts}) which is rather unusual and very strong.  The families that appear naturally in applications tend to be defined on the full parameter space $X$ and are continuous with respect to either the flat or the ${\bf F}$ topologies. The use of discretizations for families in $\mathcal Z_n(M;G)$ seems to be essential in order to prove regularity, but it can be restricted to the proofs so that a more reasonable continuous theory can be formulated. This is what we accomplish in this section.

Let $\Phi:X \rightarrow \mathcal Z_n(M^{n+1};{\bf F};G)$ be a continuous map. We let $\Pi$ be the class of all continuous maps $\Phi':X \rightarrow \mathcal Z_n(M^{n+1};{\bf F};G)$ such that $\Phi$ and $\Phi'$ are homotopic to each other in the flat topology. 

\subsection{Definition} The {\it width} of
$\Pi$ is defined by:
 $$
 {\bf L}(\Pi) = \inf_{\Phi \in \Pi}\sup_{x\in X}\{{\bf M}(\Phi(x))\}.
 $$
 If $X=I^k$, we can also restrict to maps that are zero at $\partial I^k$ and homotopies are taken relative to $\partial I^k$.
 
 \subsection{Example} Let $f:M \rightarrow \mathbb{R}$ be a Morse function with $f(M)=[0,1]$ and set $\Phi(t)=\partial \left(\{x\in M: f(x)<t\}\right)$, $t\in [0,1]$. If $M$ is non-orientable we must necessarily choose $G=\mathbb{Z}_2$, otherwise we have both possibilities $G=\mathbb{Z}$ or $G=\mathbb{Z}_2$. Then $\Phi:I\rightarrow \mathcal Z_n(M^{n+1};{\bf F};G)$ is continuous, with $\Phi(0)=\Phi(1)=0$. We denote by $\Pi_1$ the class of all maps $\Phi':I\rightarrow \mathcal Z_n(M^{n+1};{\bf F};G)$ such that $\Phi'(0)=\Phi'(1)=0$ and such that $\Phi'$ is homotopic to $\Phi$ in the flat topology relative to $\partial I=\{0,1\}$. Almgren's work \cite{almgren} (see also Section 3 of \cite{marques-neves-infinitely}) implies ${\bf L}(\Pi)>0$.

\medskip

\subsection{Definition} A sequence $\{\Phi_i\}_i\subset \Pi$ is called a {\it min-max sequence} if 
$$
{\bf L}(\Phi_i):=\sup_{x\in X}{\bf M}(\Phi_i(x))
$$
satisfies ${\bf L}(\{\Phi_i\}_{i\in \mathbb{N}}):= \limsup_{i\rightarrow \infty} {\bf L}(\Phi_i)={\bf L}(\Pi)$.

\subsection{Definition}\label{image.set} The {\it image set} of $\{\Phi_i\}_{i\in\N}$ is defined by
\begin{eqnarray*}
{\bf \Lambda}(\{\Phi_i\}_{i\in\N})&=& \{ V \in \mathcal V_n(M): \exists {\rm \, sequences\, } \{i_j\}\rightarrow \infty, x_{i_j}\in X\\
&&  \hspace{1cm} {\rm such \, that} \lim_{j\rightarrow \infty} {\bf F}(|\Phi_{i_j}(x_{i_j})|,V)=0\}.
\end{eqnarray*}

\subsection{Definition}\label{critical.set} If $\{\Phi_i\}_{i\in\N}$ is a  min-max sequence in $\Pi$ such that $L={\bf L}(\{\Phi_i\}_{i\in\N})$, the {\it critical set} 
of $\{\Phi_i\}_{i\in\N}$ is defined by
$$
{\bf C}(\{\Phi_i\}_{i\in\N})=\{V\in {\bf \Lambda}(\{\Phi_i\}): ||V||(M)=L\}.
$$

\subsection{Pull-tight} Following Pitts (\cite{pitts} p.153, see also \cite{colding-delellis, marques-neves}), we can define a continuous map 
\begin{eqnarray*}
H: I \times (\mathcal Z_n(M^{n+1};{\bf F};G) &\cap& \{T: {\bf M}(T)\leq 2{\bf L}(\Pi)\}) \\
&\rightarrow& (\mathcal Z_n(M^{n+1};{\bf F};G) \cap \{T: {\bf M}(T)\leq 2{\bf L}(\Pi)\})
\end{eqnarray*}
such that:
\begin{itemize}
\item $H(0,T)=T$ for all $T$;
\item $H(t,T)=T$ for all $t\in [0,1]$ if $|T|$ is stationary;
\item ${\bf M}(H(1,T))< {\bf M}(T)$ if $|T|$ is not stationary.
\end{itemize}

Given a min-max sequence $\{\Phi_i^*\}\subset \Pi$, we define $\Phi_i(x)=H(1,\Phi_i^*(x))$ for every $x\in X$. Then $\{\Phi_i\}\subset \Pi$ is also
a min-max sequence. It follows from the construction that ${\bf C}(\{\Phi_i\}) \subset {\bf C}(\{\Phi_i^*\})$ and that every element of 
${\bf C}(\{\Phi_i\})$ is stationary. 

\subsection{Definition} Any min-max sequence $\{\Phi_i\}\subset \Pi$ such that every element of ${\bf C}(\{\Phi_i\})$ is stationary is called {\it pulled-tight.}

\subsection{Min-max Theorem}\label{minmax.continuous.thm} {\em Suppose ${\bf L}(\Pi)>0$, and let $\{\Phi_i\}_i$ be a min-max sequence for $\Pi$. Then there exists a stationary integral varifold  $V\in {\bf C}(\{\Phi_i\})$ (hence  $||V||(M)={\bf L}(\Pi)$), with support a  closed  minimal hypersurface that is smooth embedded outside a set of dimension $n-7$.}

\subsection{Remark:} A slightly different version of this theorem was proven in Theorem 2.4 of \cite{marques-neves-cdm}. 

\begin{proof}
It follows from the pull-tight procedure that we can assume every element of ${\bf C}(\{\Phi_i\})$ is stationary. Given $\Phi_i: X \rightarrow \mathcal Z_n(M^{n+1};{\bf F};G)$, the continuity in the ${\bf F}$-metric and a compactness argument  imply there is no concentration of mass:
$$
\lim_{r\rightarrow 0}\sup_{x\in X, p\in M}||\Phi_i(x)||(B_r(p))=0,
$$
where $B_r(p)$ denotes the open geodesic ball of radius $r$ and center $p\in M$.
Theorem 3.9 of \cite{marques-neves-infinitely} (the analogous statement holds for $G=\mathbb{Z}$) implies there exist a sequence of maps
$$\phi_i^j:X(k_j^i)_0 \rightarrow \mathcal{Z}_n(M;G),$$
with $k_j^i<k_{j+1}^i$, and a
sequence of positive numbers $\{\delta_j^i\}_{j\in\N}$ converging to zero such that
\begin{itemize}
\item[(i)] $$S_i=\{\phi_i^j\}_{j\in\N}$$ is an $(X,{\bf M})$-homotopy sequence of mappings into $\mathcal{Z}_n(M;{\bf M};G)$ with ${\bf f}(\phi_i^j)<\delta_j^i$;
\item[(ii)] $$\sup\{\mathcal F(\phi_i^j(x)-\Phi_i(x)): x\in X(k_j^i)_0\}\leq \delta_j^i;$$
\item[(iii)]$$\sup\{{\bf M}(\phi_i^j(x)): x\in X(k_j^i)_0\}\leq \sup\{{\bf M}(\Phi_i(x)): x\in X\}+\delta_j^i.$$
\end{itemize}
We can also guarantee that (see Theorem 13.1 (i) of \cite{marques-neves} for the case $X=I^m$)
\begin{itemize}
\item[(iv)] for some $l_j^i\rightarrow \infty$  and every $y\in X(k_j^i)_0$,
$$
{\bf M}(\phi_i^j(y))\leq \sup \{{\bf M}(\Phi_i(x)):\alpha \in X(l_j^i), x,y\in \alpha\} + \delta_j^i.
$$
\end{itemize}

Because $x\in X \mapsto {\bf M}(\Phi_i(x))$ is continuous, we get from property (iv) above that for every $y\in X(k_j^i)_0$,
$$
{\bf M}(\phi_i^j(y)) \leq {\bf M}(\Phi_i(y))+\eta_j^i
$$
with $\eta_j^i\rightarrow 0$ as $j \rightarrow \infty$.
If we apply Lemma 4.1 of \cite{marques-neves} with $\mathcal S=\Phi_i(X)$,  property [(ii)] above implies
$$\sup\{{\bf F}(\phi_i^j(x),\Phi_i(x)): x\in X(k_j^i)_0\} \rightarrow 0$$
as $j \rightarrow \infty$.  

We can choose $j(i) \rightarrow \infty$ as $i \rightarrow \infty$ (note that necessarily $k_{j(i)}^i\rightarrow \infty$) such that  $\varphi_i=\phi_i^{j(i)}: X(k_{j(i)}^i)_0 \rightarrow \mathcal Z_n(M;G)$
satisfies:
\begin{itemize}
\item $\sup\{{\bf F}(\varphi_i(x),\Phi_i(x)): x\in X(k_{j(i)}^i)_0\} \leq a_i$ with $\lim_{i\rightarrow \infty}a_i=0$;
\item $\sup\{{\bf F}(\Phi_i(x),\Phi_i(y)): x,y \in \alpha, \alpha \in X(k_{j(i)}^i)\} \leq a_i$;
\item the finess ${\bf f}(\varphi_i)$ tends to zero as $i \rightarrow \infty$;
\item the Almgren extension $\Phi_i^{j(i)}: X \rightarrow \mathcal Z_n(M; {\bf M}; G)$ is homotopic to $\Phi_i$ in the flat topology
(by Corollary 3.12 (ii) of \cite{marques-neves-infinitely} and its analogue for $G=\mathbb{Z}$).
\end{itemize}
Following the notation of Subsection 2.12 of \cite{marques-neves-infinitely}, if $S=\{\varphi_i\}_i$ we have  ${\bf L}(S)={\bf L}(\{\Phi_i\}_{i\in \mathbb{N}})$ and 
${\bf C}(S)={\bf C}(\{\Phi_i\}_{i\in \mathbb{N}})$. In particular, every element of ${\bf C}(S)$ is stationary.

 By Theorem 2.13 of \cite{marques-neves-infinitely} (and its analogue for $G=\mathbb{Z}$), if no element $V\in {\bf C}(S)$ is $G$-almost minimizing in annuli we can find a sequence $S^*=\{\varphi_i^*\}$ of maps
 $$
 \varphi_i^*:X(k_{j(i)}^i+l_i)_0 \rightarrow \mathcal Z_n(M;G)
 $$
such that
\begin{itemize}
\item $\varphi_i$ and $\varphi_i^*$ are homotopic to each other with finesses tending to zero as $i\rightarrow \infty$;
\item ${\bf L}(S^*)< {\bf L}(S)$.
\end{itemize}

By Corollary 3.12 (i) of \cite{marques-neves-infinitely} (and its analogue for $G=\mathbb{Z}$), the Almgren extensions of $\varphi_i$, $\varphi_i^*$:
$$
\Phi_i^{j(i)}, \Phi_i^*: X \rightarrow \mathcal Z_n(M; {\bf M}; G),
$$
respectively, are homotopic to each other in the flat topology for sufficiently large $i$. Hence $\Phi_i^*$ is homotopic to $\Phi_i$ in the flat topology and so $\Phi_i^*\in \Pi$ for sufficiently large $i$. It follows from Theorem 3.10 of \cite{marques-neves-infinitely} (and its analogue for $G=\mathbb{Z}$) that 
$$
\limsup_{i\rightarrow \infty} \sup \{{\bf M}(\Phi_i^*(x): x\in X\} \leq {\bf L}(S^*) < {\bf L}(S)={\bf L}(\Pi).
$$
This is a contradiction, hence some element $V\in {\bf C}(S)$ is $G$-almost minimizing in annuli. Since ${\bf C}(S)={\bf C}(\{\Phi_i\}_{i\in \mathbb{N}})$, the theorem follows by the regularity theory of Pitts \cite[Theorem 3.13 and Section 7]{pitts} and Schoen-Simon \cite[Theorem 4]{schoen-simon} (see also Theorem 2.11 of \cite{marques-neves-infinitely}).
\end{proof}

The following remark will be useful  in Section \ref{index.bounds} where we handle the case in which the min-max hypersurfaces are one-sided.

\subsection{Remark}\label{xin.zhou} The proof actually gives more. It was observed by X. Zhou \cite{zhou} that Pitts contradiction arguments (Theorem 4.10 of \cite{pitts}) imply that some element $V$ of ${\bf C}(S)={\bf C}(\{\Phi_i\}_{i\in \mathbb{N}})$ must be $G$-almost minimizing in annuli in the following stronger sense (compared to Pitts notion): for every $p\in M$ there exists $r>0$ so that for any $0<s<r$ and any $\varepsilon>0$ there exists 
$T \in \mathcal Z_n(M;G)$ and $\delta>0$ such that $T \in \mathfrak{A}_n(B_r(p)\setminus \bar B_s(p), \varepsilon, \delta)$ (Definition 6.3 of \cite{zhou}, Definition 3.1 of \cite{pitts}) and ${\bf F}(V,|T|)< \varepsilon$. If $\Phi(x)$ is a boundary for any $\Phi\in \Pi$ and $x\in X$, we can further require that $T$ is a boundary too, i.e., there exists $U\in {\bf I}_{n+1}(M;G)$ such that $T=\partial U$. In that case we say $V$ is $G$-almost minimizing of {\it boundary type} in annuli.

\section{Unstable variations and the Morse index}

If $\Sigma$ is a smooth embedded closed minimal hypersurface, its {\it Morse index}  is the maximal dimension of a subspace
of normal variations restricted to which the second variation of area quadratic form $Q(\cdot, \cdot)$ is negative definite. Here we have to consider
the general case of stationary varifolds.

\subsection{Definition}\label{unstable} Let $\Sigma$ be a stationary varifold in $ \mathcal{V}_n(M)$ and $\varepsilon\geq 0$. We say   that $\Sigma$ is {\em $k$-unstable in an $\varepsilon$-neighborhood} if there exist $0<c_0<1$ and  a smooth family $\{F_v\}_{v\in \overline B^k}\subset \text{Diff}(M)$ with $F_0={\rm Id}$, $F_{-v}=F_v^{-1}$ for all $v\in\overline B^k$  such that, for any $V\in {\overline{\bf B}^{\bf F}_{2\varepsilon}(\Sigma)},$ the smooth function
$$A^V:\overline B^k\rightarrow [0,\infty),\quad\quad A^V(v)=||(F_v)_{\#}V||(M)$$ satisfies:
\begin{itemize}
\item $A^V$ has a unique maximum at $m(V)\in B^k_{c_0/\sqrt{10}}(0)$;
\item 
$-\frac{1}{c_0}\,{\rm Id}\leq {D^2 A^{V}}(u)\leq -c_0\,{\rm Id}\quad\mbox{for all }u\in \overline B^k.$
\end{itemize}
Here $(F_v)_{\#}$ denotes the push-forward operation. 
Also, because $\Sigma$ is stationary, necessarily $m(\Sigma)=0$.  


\medskip

The map $V\mapsto m(V)$ is continuous on $\overline{\bf B}^{\bf F}_{2\varepsilon}(\Sigma)$ and 
\begin{equation}\label{area.inequality}
 A^V(m(V))-\frac{1}{2c_0}{|u-m(V)|^2}\leq A^V(u)\leq A^V(m(V))-\frac{c_0}{2}{|u-m(V)|^2}
 \end{equation}
for all $u\in \overline B^k$. 

If $V_i$ tends to $V$ in the ${\bf F}$-topology then $A^{V_i}$ tends to $A^V$ in the smooth topology \cite[Section 2.3]{pitts}. Thus if $\Sigma$ is stationary and $k$-unstable in a $0$-neighborhood then it is $k$-unstable in an $\varepsilon$-neighborhood for some $\varepsilon>0$. 

\subsection{Definition}
We say that $\Sigma$ is {\em $k$-unstable} if it is stationary and $k$-unstable in an 
$\varepsilon$-neighborhood for some $\varepsilon>0$. 
\medskip

\subsection{Proposition} {\em If  $\Sigma$ is smooth embedded closed and minimal, then $\Sigma$ is $k$-unstable if and only if   ${\rm index}(\Sigma)\geq k$.} 

\begin{proof}
Denote by $Q(\cdot,\cdot)$ the second variation of area quadratic form, acting on normal vector fields along $\Sigma$. Suppose ${\rm index}(\Sigma)\geq k$. Then there exist normal vector fields $X_1,\dots,X_k$ along $\Sigma$ such that 
$$
Q(\sum_ia_iX_i,\sum_ia_iX_i)<0
$$ 
for any $(a_i)\neq 0\in \mathbb{R}^k$. These vector fields can be extended to ambient vector fields $X_1,\dots,X_k$  in $M$. Denote by $\phi_t^X$ the flow of $X$ and, for some $\delta$ to be chosen later, define
$$
F_v(x)=F(v,x)=\phi^{\sum_iv_iX_i}_{\delta}(x)
$$
for $v \in \overline B^k$ and $x\in M$. Then $F_0={\rm Id}$, $F_{-v}=F_v^{-1}$ for all $v\in\overline B^k$, $DA^\Sigma(0)=0$.
Moreover, since
$$
F(tv,x)=\phi^{t\sum_iv_iX_i}_{\delta}(x)=\phi^{\delta \sum_iv_iX_i}_{ t}(x),
$$
we have 
$$
\frac{\partial}{\partial t}_{|t=0}F(tv,x)=\delta \sum_iv_iX_i(x).
$$
Hence $D^2A^\Sigma(0)<0$. By choosing $\delta$ sufficiently small, we conclude that $\Sigma$ is $k$-unstable. 

Now suppose $\Sigma$ is $k$-unstable, and let $\{F_v\}_{v\in \overline B^k}$ be the corresponding family of diffeomorphisms.  Define 
$Y_i(x)=\frac{\partial}{\partial t}_{|t=0}F(te_i,x)=\frac{\partial F}{\partial v_i}(0,x)$. Then 
$$
0> \frac{d^2}{dt^2}_{|t=0}{\rm area}(F_{tv}(\Sigma))=Q (\sum_i v_iY_i^\perp,\sum_i v_iY_i^\perp),
$$
where $X_i=Y_i^\perp$ is the normal component of $Y_i$  along $\Sigma$. This means $\{X_1,\dots,X_k\}$ is linearly independent and spans a subspace of normal vector fields restricted to which $Q$ is negative definite. Therefore ${\rm index}(\Sigma)\geq k$.

\end{proof}

In what follows we assume that $\Sigma$ is a stationary $n$-varifold that is $k$-unstable in an $\varepsilon$-neighborhood, $\varepsilon>0$,  and derive two basic lemmas. The family $\{F_v\}_{v\in \overline B^k}$ is as in Definition \ref{unstable}.

\subsection{Lemma}\label{away.lemma}{\em There exist $\bar\eta=\bar \eta(\varepsilon, \Sigma,  \{F_v\})>0$, so that if $V\in \mathcal V_n(M)$ with ${\bf F}(V,\Sigma)\geq \varepsilon$ satisfies
$$||(F_v)_{\#}V||(M)\leq ||V||(M)+\bar \eta$$
for some $v \in \overline{B}^k$, then ${\bf F}((F_v)_{\#}V,\Sigma)\geq 2\bar\eta.$
}

\begin{proof}

We argue by contradiction and assume that there is a sequence $V_i\in \mathcal V_n(M)$  with ${\bf F}(V_i,\Sigma)\geq \varepsilon$, $$||(F_{v_i})_{\#}V_i||(M)\leq ||V_i||(M)+\frac 1 i$$ for some $v_i \in \overline B^k$, and $(F_{v_i})_{\#}V_i$ tending to $\Sigma$ in the ${\bf F}$-topology. After passing to a subsequence we can assume that $v_i$ tends to $v\in \overline B^k$. Hence $V_i$ converges to $(F_{-v})_{\#}\Sigma$ and  $||\Sigma||(M)\leq ||(F_{-v})_{\#}\Sigma||(M)$, which forces $v=0$. Thus $V_i$ tends to $\Sigma$, which is impossible because  ${\bf F}(V_i,\Sigma)\geq \varepsilon$ for all $i\in\N$.
\end{proof}

For each $V\in \overline{\bf B}^{\bf F}_{2\varepsilon}(\Sigma)$, consider the one-parameter flow   $\{\phi^V(\cdot,t)\}_{t\geq 0}\subset \mathrm{Diff}(\overline B^k)$  generated by the vector field $$u\mapsto -(1-|u|^2)\nabla A^V(u), \quad  u\in \overline B^k.$$
With $u\in \overline B^k$ fixed we have that $t\mapsto A^V(\phi^V(u,t))$ is non-increasing.

\subsection{Lemma}\label{T.lemma}{\em For all $\delta< 1/4$ there is $T=T(\delta,\varepsilon, \Sigma, \{F_v\}, c_0)\geq 0$ so that for any $V\in \overline{\bf B}^{\bf F}_{2\varepsilon}(\Sigma)$ and $v\in \overline B^k$ with $|v-m(V)|\geq \delta$ we have
$$ A^V(\phi^V(v,T))<  A^V(0)-\frac{c_0}{10}\quad\mbox{and}\quad |\phi^V(v,T)|>\frac{c_0}{4}.$$
}

\begin{proof}
Because $m(V)\in B^k_{c_0/\sqrt {10}}(0)$,  we have by inequality (\ref{area.inequality}) that
$$\sup_ {|u|\leq 1} A^V(u)=A^V(m(V))\leq A^V(0)+\frac{c_0}{20}.$$
So, to prove the first inequality in Lemma \ref{T.lemma} it suffices to show  the existence of $T$ so that
\begin{equation}\label{T.lemma.inequality}
 V\in \overline{\bf B}^{\bf F}_{2\varepsilon}(\Sigma) , |v-m(V)|\geq \delta \implies A^V(\phi^V(v,T))<\sup_ {|u|\leq 1} A^V(u)-\frac{c_0}{5}.
\end{equation}

We argue by contradiction and assume that  for all $i\in\N$ there exist $v_i\in \overline B^k$ and $V_i\in \overline{\bf B}^{\bf F}_{2\varepsilon}(\Sigma)$ so that $|v_i-m(V_i)| \geq\delta$ but
$$ A^{V_i}(\phi^{V_i}(v_i,i))\geq  \sup_ {|u|\leq 1} A^{V_i}(u)-\frac{c_0}{5}.$$
This inequality implies 
$$ A^{V_i}(\phi^{V_i}(v_i,t))\geq  \sup_ {|u|\leq 1} A^{V_i}(u)-\frac{c_0}{5}\quad\text{for all }0\leq t\leq i.$$

Since both $\overline{\bf B}^{\bf F}_{2\varepsilon}(\Sigma)$ and $\overline B^k$ are compact, we obtain subsequential limits $V\in \overline{\bf B}^{\bf F}_{2\varepsilon}(\Sigma)$ and $v\in \overline B^k\setminus B^k_{\delta}(m(V))$ so that
$$ A^{V}(\phi^{V}(v,t))\geq \sup_ {|u|\leq 1} A^{V}(u)-\frac{c_0}{5} \quad\text{for all }t\geq 0.$$

Because $|v-m(V)|\geq \delta$, we have $\lim_{t\to\infty}|\phi^V(v,t)|=1$ and thus we deduce from the inequality above that
$$\sup_ {|u|=1} A^V(u)\geq \sup_ {|u|\leq 1} A^{V}(u)-\frac{c_0}{5}.$$
On the other hand, $|u-m(V)|>2/3$ for all $u\in \overline B^k$ with $|u|=1$. Hence, by inequality (\ref{area.inequality}),
$$\sup_ {|u|=1} A^V(u)\leq\sup_ {|u|\leq 1} A^{V}(u)-\frac {c_0}{2}\left(\frac{2}{3}\right)^2<\sup_ {|u|\leq 1} A^{V}(u)-\frac{c_0}{5},$$
which gives us the desired contradiction.

We now show  the second inequality of Lemma \ref{T.lemma}. If $w\in \overline B^k$ has $|w|\leq c_0/4$, then $|w-m(V)|\leq c_0(1/4+1/\sqrt{10})\leq c_0(1/4+1/3)$ and so, using inequality (\ref{area.inequality}) again,
$$A^V(w)\geq \sup_ {|u|\leq 1} A^{V}(u)-\frac{c_0}{2}\left(\frac{1}{4}+\frac{1}{3}\right)^2>\sup_ {|u|\leq 1} A^{V}(u)-\frac{c_0}{5}.$$
Combining this inequality with \eqref{T.lemma.inequality} we obtain the desired statement.
\end{proof}

\section{Deformation theorems}\label{deformations}

Let  $\{\Phi_i\}_{i\in\N}$ be a  sequence of continuous maps from $X$ into $\mathcal Z_n(M;{\bf F};G)$. We set 
$$L={\bf L}(\{\Phi_i\}_{i\in\N}):=\limsup_{i \rightarrow \infty} \sup_{x\in X} {\bf M}(\Phi_i(x)).$$ If $X=I^k$, we can also consider maps that are zero at $\partial I^k$ with homotopies  relative to $\partial I^k$.

In this section we will prove some deformation theorems that will be fundamental in the proof of the index estimates. The basic idea is to produce
homotopic families $\{\Psi_i\}_{i\in\N}$ to $\{\Phi_i\}_{i\in\N}$ that do not increase the area and avoid minimal hypersurfaces with large index or multiplicity.

\subsection{Deformation Theorem A}\label{induction.theorem}{\em Suppose that 
\begin{enumerate}
\item[(a)] $\Sigma\in \mathcal V_n(M)$ is stationary and $(k+1)$-unstable;
\item[(b)] $K\subset \mathcal V_n(M)$  is  a compact set  so that $\Sigma\notin K$ and $|\Phi_i|(X)\cap K=\emptyset$ for all $i\geq i_0$;
\item[(c)] $||\Sigma||(M)= L$.
\end{enumerate}
There exist  $\bar \varepsilon>0$, $j_0\in\N$, and another sequence $\{\Psi_i\}_{i\in\N}$ so that 
\begin{itemize}
\item[(i)] $\Psi_i$ is homotopic to $\Phi_i$ in the ${\bf F}$-topology for all $i\in\N$;
\item[(ii)] $ {\bf L}(\{\Psi_i\}_{i\in\N})\leq L$;
\item[(iii)] for all $i\geq j_0$, $|\Psi_i|(X)\cap ({\overline{\bf B}}^{\bf F}_{\bar \varepsilon}(\Sigma)\cup K)=\emptyset$.
\end{itemize}
}

\begin{proof}
Set $d=\min\{{\bf F}(\Sigma, Z):Z\in K\}>0.$

By assumption, $\Sigma$ satisfies the conditions of Definition \ref{unstable} with $(k+1)$ in place of $k$  for some $\varepsilon>0$. Let $\{F_v\}_{v\in \overline B^{k+1}}$ and $c_0$  be the corresponding family and constant.
 Without loss of generality (by changing $\varepsilon$, $\{F_v\}$, $c_0$) we can assume that
 \begin{equation}\label{K.disjoint}
 \min\{{\bf F}((F_{v})_{\#}V, Z):Z\in K, v\in \overline B^{k+1}\}>d/2 \quad\mbox{for all }V\in \overline{\bf B}^{\bf F}_{2\varepsilon}(\Sigma).
 \end{equation}
 
 Let $X(k_i)$ be a sufficiently fine subdivision of $X$ (see Section 2 of \cite{marques-neves-infinitely}) so that ${\bf F}(|\Phi_i(x)|,|\Phi_i(y)|)<\delta_i$ for any $x,y$ belonging to the same cell in $X(k_i)$, with $\delta_i=\min \{2^{-(i+k+2)}, \varepsilon/4\}$. We can also assume $$|m(|\Phi_i(x)|)-m(|\Phi_i(y)|)|<\delta_i$$
for any $x,y$ with ${\bf F}(|\Phi_i(x)|,\Sigma)\leq 2\varepsilon$, ${\bf F}(|\Phi_i(y)|,\Sigma)\leq 2\varepsilon$, and belonging to the same cell in $X(k_i)$.  For $\eta>0$, let $U_{i,\eta}$ be the union of all cells $\sigma \in X(k_i)$ so that ${\bf F}(|\Phi_i(x)|,\Sigma)<\eta$ for all $x\in \sigma$. Then $U_{i,\eta}$ is a subcomplex of $X(k_i)$. If $\beta \notin U_{i,\eta}$, then ${\bf F}(|\Phi_i(x')|,\Sigma)\geq \eta$ for some $x'\in \beta$. Consequently,  ${\bf F}(|\Phi_i(x)|,\Sigma)\geq \eta-\delta_i$ for all $x\in \beta$. For each $i\in\N$ and $x\in  U_{i,2\varepsilon}$ we use the notation $A_i^x=A^{|\Phi_i(x)|}$, $m_i(x)=m(|\Phi_i(x)|)$ and $\phi_i^x=\phi^{|\Phi_i(x)|}$.

 The function $m_i:U_{i,2\varepsilon}\rightarrow \overline B^{k+1}$  is continuous.  We would like to define a continuous homotopy
$$\hat H_i:U_{i,2\varepsilon}\times  [0,1]\rightarrow B^{k+1}_{1/2^i}(0) \quad\text{so that } \hat H_i(x,0)=0\text{ for all }x \in U_{i,2\varepsilon}$$
and
\begin{equation}\label{delta.definition}
\inf_{x\in U_{i,2\varepsilon}}|m_i(x)-\hat H_i(x,1)|\geq \eta_i>0
\end{equation}
for some $\eta_i>0$. Hypothesis (a) is crucial for the construction of $\hat H_i$ and, roughly speaking, the motivation is the following: the subspaces $A_i=\{(x,0):x\in U_{i,2\varepsilon}\}$ and $B_i=\{(x,m_i(x)):x\in U_{i,2\varepsilon}\}$  have both dimension $k$ and are contained in a space of dimension $2k+1$. Hence one should be able to perturb $A_i$ slightly so that it does not intersect $B_i$ and thus show \eqref{delta.definition}. We now proceed to describe the construction of $\hat H_i$ carefully.

Because $B^{k+1}_{1/2^i}(0)$ is convex it is enough to define $a_i(x)=\hat H_i(x,1)$. We will write $a_i(x)=m_i(x)+b_i(x)$ and
define $b_i(x)$ for $x \in U_{i,2\varepsilon}$. Hence we need $b_i(x) \neq 0$ and $b_i(x) \in B^{k+1}_{1/2^i}(-m_i(x))$ for all $x \in U_{i,2\varepsilon}$.

We will define the function $b_i$ inductively on the skeletons of $U_{i,2\varepsilon}$. Let $U_{i,2\varepsilon}^{(j)}$ denote the $j$-dimensional skeleton, and recall that $X$ is $k$-dimensional. For every $x\in U_{i,2\varepsilon}^{(0)}$ it is easy to choose $b_i(x)\neq 0$ and $b_i(x) \in B^{k+1}_{2^{-(i+k)}}(-m_i(x))$. Suppose we have defined $b_i$ on $U_{i,2\varepsilon}^{(j)}$ so that $b_i(x)\neq 0$ and $$b_i(x) \in B^{k+1}_{2^{-(i+k-j)}}(-m_i(x)) {\rm \, \, for \, \, all\, } x\in U_{i,2\varepsilon}^{(j)}. $$
Let $\sigma \in U_{i,2\varepsilon}^{(j+1)}$, and let $x_\sigma$ be its center. Then $b_i(x) \in B^{k+1}_{2^{-(i+k-j)}+\delta_i}(-m_i(x_\sigma))$ for every $x \in \partial \sigma$. Because ${\rm dim}(\partial \sigma)\leq j\leq (k-1)$ and $\pi_j(B^{k+1}\setminus\{0\})=0$, we can find a continuous extension
of $b_i$ to $\sigma$ so that $b_i(x)\neq 0$ and $b_i(x) \in B^{k+1}_{2^{-(i+k-j)}+\delta_i}(-m_i(x_\sigma))$ for every $x \in  \sigma$. Then
$b_i(x) \in B^{k+1}_{2^{-(i+k-j)}+2\delta_i}(-m_i(x))$ for every $x \in  \sigma$. But $2^{-(i+k-j)}+2\delta_i< 2^{-(i+k-j-1)}$, completing the induction step.

Let $c:[0,\infty]\rightarrow [0,1]$ be a cutoff function which is non-increasing, one in a neighborhood of $[0,3\varepsilon/2]$, and zero in a neighborhood of $[7\varepsilon/4, +\infty)$. Note that if $y \notin U_{i,2\varepsilon}$, then ${\bf F}(|\Phi_i(y)|,\Sigma) \geq 2\varepsilon-\delta_i \geq 7\varepsilon/4$. Hence $$c({\bf F}(|\Phi_i(y)|,\Sigma))=0  {\rm \, \, for \, \, all\, } y \notin U_{i,2\varepsilon}.$$

Consider
$H_i:X\times[0,1]\rightarrow B^{k+1}_{2^{-i}}(0)$ given by
$$H_i(x,t)=\hat H_i(x,c({\bf F}(|\Phi_i(x)|,\Sigma))t)\text{ if }x\in U_{i,2\varepsilon}$$
and $$H_i(x,t)=0\text {\, \, if }x\in X\setminus U_{i,2\varepsilon}.$$
The map $H_i$ is continuous.

With $\eta_i$ given by \eqref{delta.definition}, consider $T_i=T(\eta_i,\varepsilon, \Sigma, \{F_v\}, c_0)\geq 0$ given by Lemma \ref{T.lemma}. We now set
$$D_i:X\rightarrow \overline B^{k+1}, \quad D_i(x)=\phi_i^x(H_i(x,1), c({\bf F}(|\Phi_i(x)|,\Sigma))T_i)\text{ if }x\in U_{i,2\varepsilon}$$ and
$$D_i(x)=0\text {\, \, if }x\in X\setminus U_{i,2\varepsilon}.$$
The map $D_i$ is continuous.
Define $$\Psi_i:X\rightarrow   \mathcal Z_n(M;{\bf F};G), \quad \Psi_i(x)=(F_{D_i(x)})_{\#}(\Phi_i(x)).$$ In particular,
$$\Psi_i(x)=\Phi_i(x) \text {\, \, if }x\in X\setminus U_{i,2\varepsilon}.$$ 

 The map $D_i$ is homotopic to the zero map in $\overline B^{k+1}$, so $\Psi_i$ is homotopic to $\Phi_i$ in the ${\bf F}$-topology for all $i\in\N$.

\medskip

\noindent{\bf Claim 1:} ${\bf L}(\{\Psi_i\}_{i\in\N})\leq L$.\\
 
From the mass non-increasing property of $\phi^x_i$ we have that for all $x\in X$
$$
||\Psi_i(x)||(M)\leq ||(F_{H_i(x,1)})_{\#}(\Phi_i(x))||(M).
$$
Note that $||F_v-\text{Id}||_{C^2}$ tends to zero uniformly as  $v\in B^{k+1}$ tends to zero. Thus, using the fact that $H_i(x,1)\subset B_{1/2^i}^{k+1}(0)$ for all $x\in X$, we have that
\begin{equation}\label{uniform.mass}
\lim_{i\to\infty}\sup_{x\in X} \left|\,||\Phi_i(x)||(M)-||(F_{H_i(x,1)})_{\#}(\Phi_i(x))||(M)\,\right|=0
\end{equation}
and  this finishes the proof.
\medskip

\noindent{\bf Claim 2: }{\em There exists $\bar\varepsilon>0$ so that, 
for all  sufficiently large $i$, $|\Psi_i|(X)\cap \overline {\bf B}^{\bf F}_{\bar \varepsilon}(\Sigma)=\emptyset$.}\\

There are  three cases to consider.
If $x\in X\setminus U_{i,2\varepsilon}$, then $\Psi_i(x)=\Phi_i(x)$ and so ${\bf F}(|\Psi_i(x)|,\Sigma)\geq 7\varepsilon/4$. 

 If $x\in U_{i,2\varepsilon}\setminus U_{i,5\varepsilon/4}$, then ${\bf F}(|\Phi_i(x)|,\Sigma)\geq \varepsilon$. The mass non-increasing property of $\phi^x_i$ implies
$$||\Psi_i(x)||(M)=||(F_{D_i(x)})_{\#}(\Phi_i(x))||(M) \leq ||(F_{H_i(x,1)})_{\#}(\Phi_i(x))||(M).$$
From \eqref{uniform.mass} we have that for all $i$ sufficiently large
$$ ||(F_{H_i(y,1)})_{\#}(\Phi_i(y))||(M)\leq ||\Phi_i(y)||(M)+\bar\eta\quad\text{for all }y\in X,$$
where $\bar\eta=\bar \eta(\varepsilon, \Sigma,  \{F_v\})>0$ is given by Lemma \ref{away.lemma}. Combining the two inequalities above with Lemma \ref{away.lemma} applied to $V=|\Phi_i(x)|$ and $v=D_i(x)$ we have that ${\bf F}(|\Psi_i(x)|,\Sigma)\geq 2\bar\eta$.

The last case to consider is when  $x\in U_{i,5\varepsilon/4}$.  Then $c({\bf F}(|\Phi_i(x)|,\Sigma))=1$. Note that there exists $\bar \gamma=\bar \gamma(\Sigma,  c_0)>0$  so that 
\begin{equation}\label{claim2.inequality}
||V||(M)\leq ||\Sigma||(M)-\frac{c_0}{20}\implies {\bf F}(V,\Sigma)\geq 2\bar\gamma.
\end{equation}
 Choose $i$ large enough so  that 
$$\sup_{x\in X}||\Phi_i(x)||(M)\leq ||\Sigma||(M)+\frac{c_0}{20}.$$
From Lemma \ref{T.lemma} (with $\delta=\eta_i$, $V=|\Phi_i(x)|$, $v=H_i(x,1)$) we have
$$||\Psi_i(x)||(M)=A_i^{x}(\phi_i^x(H_i(x,1), T_i))< A^x_i(0)-\frac{c_0}{10} =||\Phi_i(x)||(M)- \frac{c_0}{10}$$
and so 
$$||\Psi_i(x)||(M)\leq  ||\Sigma||(M)-\frac{c_0}{20}.$$
From \eqref{claim2.inequality} we obtain that ${\bf F}(|\Psi_i(x)|,\Sigma)\geq 2\bar\gamma$. This ends the proof of Claim 2.

\medskip

\noindent{\bf Claim 3: }{\em For all $i$, $|\Psi_i|(X)\cap K=\emptyset$.}\\

If $x\in X\setminus U_{i,2\varepsilon}$, then $|\Psi_i(x)|=|\Phi_i(x)|\notin K$. If $x\in U_{i,2\varepsilon}$, then ${\bf F}(|\Phi_i(x)|,\Sigma)\leq 2\varepsilon$ and we have 
  from \eqref{K.disjoint} that $|\Psi_i(x)|\notin K$. We have proved the theorem.
\medskip

\end{proof}

\subsection{Remark.} In Deformation Theorem A, if there is a subcomplex $Y\subset X$ such that $(\Phi_i)_{|Y} \equiv 0$ then $(\Psi_i)_{|Y} \equiv 0$ and the homotopy is relative to $Y$.

\medskip

Throughout the rest of the section the space $X$ can be either $S^1$ or $[0,1]$, and as before $G$ is either $\Z$ or $\Z_2$. Let $\{\Phi_i\}_{i\in\N}$ be a sequence of continuous maps from $X$ into $\mathcal Z_n(M;{\bf F};G)$ and  let $L={\bf L}(\{\Phi_i\}_{i\in\N})$.  If $X=[0,1]$, assume further that $\Phi_i(0)=\Phi_i(1)=0$ and that all homotopies are relative to $\partial [0,1]=\{0,1\}$.

\subsection{Deformation Theorem B}\label{second.index.theorem}{\em Suppose that $\Sigma\in \mathcal V_n(M)$ is a stationary integral varifold such that:
\begin{enumerate}
\item[(a)]  the support of $\Sigma$ is a 
smooth embedded closed minimal hypersurface with a connected component $S$ of Morse index  one;
\item[(b)] for any sufficiently small  tubular neighborhood $\Omega\subset M$ of $S$  we have  $||\Sigma||(\Omega)=m||S||(\Omega)$ for some integer $m\geq 2$;
\item[(c)] there exists a compact set $K\subset \mathcal V_n(M)$ so that  $\Sigma\notin K$ and $|\Phi_i|(X)\cap K=\emptyset$ for all $i\geq i_0$;
\item[(d)]$||\Sigma||(M)=L$.
\end{enumerate}
Then there exist $\xi>0$, $j_0\in\N$ so that  for all $i\geq j_0$ one can find   $\Psi_i: X\rightarrow \mathcal Z_n(M;{\bf F};G)$ such that
\begin{itemize}
\item[(i)] $\Psi_i$ is homotopic to $\Phi_i$  in the flat topology;
\item[(ii)] ${\bf L}(\{\Psi_i\}_{i\in\N})\leq L$;
\item[(iii)]  $|\Psi_i|(X)\cap {\overline{\bf B}}^{\bf F}_{\xi}(\Sigma)=\emptyset$;
\item[(iv)] $|\Psi_i|(X)\cap K\cap \{V: ||V||(M)\geq L-\xi\}=\emptyset.$
\end{itemize}

}

\begin{proof}
Set $d=\min\{{\bf F}(\Sigma, Z):Z\in K\}>0.$ Let $S_1,\dots,S_P$ be the connected components of $\Sigma$, with $S_1=S$. Then, as a stationary integral varifold, $$\Sigma = m_1\cdot |S_1| + \cdots + m_P \cdot |S_P|$$ with $m_p\in \mathbb{N}$ and $m_1\geq 2$.

Let $\mathcal S_0$ be the finite set of all $T=\sum_p a_p \cdot S_p \in \mathcal Z_n(M;G)$ with $a_p\in G$ and $|a_p|\leq m_p$. If $G=\mathbb{Z}$ and $S_p$ is non-orientable, we set $a_p=0$.
We choose a constant $\delta>0$
so that for all $T\in \mathcal S_0$, the balls ${\bf B}^{\mathcal F}_{2\delta}(T)$ are mutually disjoint and any two paths  in ${\bf B}^{\mathcal F}_{2\delta}(T)$  with common endpoints are homotopic (not
necessarily in ${\bf B}^{\mathcal F}_{2\delta}(T)$) to each other relative to  the boundary in the flat topology.
The existence of such constant follows from the work of Almgren \cite{almgren} (see also Proposition 3.5 or Appendix A of \cite{marques-neves-infinitely}). 

The fact that $S$ has  index one implies that $\Sigma$ is $1$-unstable in an $\varepsilon$-neighborhood  for some  $\varepsilon>0$ (Definition \ref{unstable}). Let $\{F_v\}_{v\in [-1,1]}$ and $c_0>0$
be the corresponding family and constant. Note that we can choose the family of diffeomorphisms $\{F_v\}$ such that they are all equal to the identity map outside some small tubular neighborhood of $S$. We 
let $J=[-1,1]$  and denote
$$\mathcal S_1=\{(F_{v})_{\#}T: v\in J, T \in \mathcal S_0\} \subset \mathcal Z_n(M;G).$$

By changing $\varepsilon>0$, $\{F_v\}$ and $c_0$ if necessary, we can assume that 
\begin{itemize}
\item[($h_0$)] every  $T\in \mathcal S_0$ with $|T| \neq \Sigma$  satisfies $||(F_{v})_{\#}(|T|)||(M)\leq ||\Sigma||(M)-c_0$ for all $v\in J$;
\item[($h_1$)] 
 for all $T\in \mathcal S_1$ and $v\in J$, $(F_{v})_{\#}T\in {\bf B}^{\mathcal F}_{\delta/(2m_1)}(T)$;
\item[($h_2$)] for all $V\in \overline {\bf B}^{\bf F}_{2\varepsilon}(\Sigma)$
 \begin{equation*} \label{K.disjoint.theorem2}
 \min\{{\bf F}((F_{v})_{\#}V, Z):Z\in K, v\in J\}>d/2;
 \end{equation*}
 \item[($h_3$)] for all $V\in \overline {\bf B}^{\bf F}_{2\varepsilon}(\Sigma)$, $||V||(M) \leq ||\Sigma||(M)+c_0/20.$
\end{itemize}

We need to prove some lemmas. 

\subsection{Lemma}\label{second.index.thm.lemma1}{\em  There exists $\varepsilon_1>0$ so that if $V\in \mathcal V_n(M)$ is such that
${\bf F}(V,U)\leq \varepsilon_1$ for some $U\in \mathcal S_1$ then $||V||(M)\leq ||U||(M)+c_0/20.$
}
\begin{proof}
The  result is a direct consequence of the fact that $|\mathcal S_1|$ is compact in $\mathcal V_n(M)$.
\end{proof}

\subsection{Lemma}\label{second.index.thm.lemma2}{\em  For all $0<\eta < 2\delta$ there exists $\varepsilon_2=\varepsilon_2(\eta,\Sigma)>0$ so that
 for any continuous map in the flat topology $\Psi:[0,1]\rightarrow \mathcal Z_n(M;G)$ with $|\Psi|([0,1])\subset {\bf B}^{\bf F}_{\varepsilon_2}(\Sigma)$ 
there exists $T \in \mathcal S_0$ so that $\Psi([0,1])\subset {\bf B}^{\mathcal F}_{\eta}(T)$.
}
\begin{proof}
First we prove that there exists $\varepsilon_2=\varepsilon_2(\eta,\Sigma)>0$ such that if $T'\in \mathcal Z_n(M;G)$ satisfies ${\bf F}(|T'|,\Sigma)<\varepsilon_2$, then there exists $T \in \mathcal S_0$ with $T'\in {\bf B}^{\mathcal F}_{\eta}(T)$.  

Suppose this is not true. Then for some $\eta>0$ we can find a sequence $T_i'\in \mathcal Z_n(M;G)$ satisfying ${\bf F}(|T_i'|,\Sigma)<1/i$ and $\mathcal F(T_i',T)\geq \eta$ for every $T\in \mathcal S_0$. We have ${\bf M}(T_i')\leq C$ for all $i$ and some $C>0$, hence by compactness a subsequence $T_j'$ converges to $\hat T\in \mathcal Z_n(M;G)$ in the flat topology. The convergence of $|T_j'|$ to $\Sigma$ in the ${\bf F}$-metric implies 
$||T_j'||(M\setminus \rm{supp}(\Sigma))\rightarrow 0$ and for appropriate tubular neighborhoods $U_p$ of $S_p$, 
$||T_j'||(U_p)\rightarrow m_p\cdot |S_p|$. Then $\rm{supp}(\hat T)\subset \rm{supp}(\Sigma)$ and by the Constancy Theorem $\hat{T}=\sum_p a_p\cdot S_p$ for $a_p\in G$, where we set $a_p=0$ if $S_p$ is unorientable and $G=\mathbb{Z}$.  Lower semicontinuity of mass implies $|a_p|\leq m_p$ and thus $\hat T \in \mathcal S_0$. Contradiction because we were assuming $\mathcal F(T_i',T)\geq \eta$ for every $T\in \mathcal S_0$.

If  $\Psi:[0,1]\rightarrow \mathcal Z_n(M;G)$ is a continuous map in the flat topology with $|\Psi|([0,1])\subset {\bf B}^{\bf F}_{\varepsilon_2}(\Sigma)$, then by what we just proved $\{\Psi^{-1}(B_\eta^{\mathcal F}(T))\}_{T\in \mathcal S_0}$ is a family of disjoint open sets covering $[0,1]$. By connectedness, there must exist a single $T \in \mathcal S_0$ so that $\Psi([0,1])\subset {\bf B}^{\mathcal F}_{\eta}(T)$.

\end{proof}

Finally, there is a constant $C\geq 1$ so that for all $S,T\in \mathcal Z_n(M;G)$ we have
$$\sup_{u\in J}\mathcal F((F_u)_{\#}S,(F_u)_{\#}T)\leq C \mathcal F(S,T).$$

Fix $\tilde\delta<\{\delta/2,c_0/40\}$. Consider $\eta=\eta(\mathcal S_1,\tilde \delta, L)<\tilde\delta$ given by Proposition \ref{f.close},  $\varepsilon_1>0$ given by Lemma \ref{second.index.thm.lemma1},  $\varepsilon_2=\varepsilon_2(\eta/C,\Sigma)>0$ given by Lemma \ref{second.index.thm.lemma2}, and choose $ 0<\bar \varepsilon \leq \min \{\varepsilon, \varepsilon_1/2,\varepsilon_2/2\}$.

If $|\Phi_i|(X)\subset {\bf B}^{\bf F}_{\varepsilon_2}(\Sigma)$,  Lemma \ref{second.index.thm.lemma2} and the choice of $\delta$ imply in this case $\Phi_i$ is homotopic to a constant path in the flat topology. Then Deformation Theorem B automatically holds by choosing $\Psi_i \equiv 0$.  Likewise, if $|\Phi_i|(X)\cap  {\bf B}^{\bf F}_{\bar\varepsilon}(\Sigma)=\emptyset$ we can just take $\Psi_i=\Phi_i$ and  so we assume that this case does not occur as well.

Elementary considerations allow us to find a closed set $V_{i,\bar \varepsilon}\subset X$ with finitely many connected components which are homeomorphic to closed intervals, $\partial V_{i,\bar \varepsilon}$ a finite union of of points, and such that
\begin{equation}\label{distance.boundary}
|\Phi_i|(V_{i,\bar \varepsilon})\subset {\bf B}^{\bf F}_{2\bar\varepsilon}(\Sigma)\quad\mbox{and}\quad |\Phi_i|(X \setminus V_{i,\bar \varepsilon}) \subset \mathcal V_n(M) \setminus  {\bf B}^{\bf F}_{\bar\varepsilon}(\Sigma).
\end{equation}
Consequently, $ |\Phi_i|(\partial V_{i,\bar \varepsilon}) \subset {\bf B}^{\bf F}_{2\bar\varepsilon}(\Sigma)\setminus  {\bf B}^{\bf F}_{\bar\varepsilon}(\Sigma).$ It also follows that ${\bf M}(\Phi_i(y)) \leq {\bf M}(\Sigma) + c_0/20$ for every $y\in V_{i,\bar \varepsilon}$  by property $(h_3)$ above.

Choose a continuous map $$ H^1_i:\partial V_{i,\bar \varepsilon}\times  [0,1]\rightarrow B^1_{1/i}(0) \quad\text{so that } H^1_i(x,0)=0\text{ for all }x \in \partial V_{i,\bar \varepsilon}$$
and
\begin{equation}\label{delta.definition2}
\inf_{x\in \partial V_{i,\bar \varepsilon}}|m_i(x)- H^1_i(x,1)|\geq \eta_i>0
\end{equation}
for some $\eta_i>0$. Consider $$D_i:\partial V_{i,\bar \varepsilon}\times  [0,1]\rightarrow J, \quad D_i(x,t)=\phi_i^x(H^1_i(x,1),tT_i),$$
where $T_i=T(\eta_i,\bar \varepsilon, \Sigma, \{F_v\},c_0)\geq 0$ is given by Lemma \ref{T.lemma}. We have $|D_i(x,1)|\geq c_0/4$ for $x\in \partial V_{i,\bar \varepsilon}$.

Let $\tilde V_{i,\bar \varepsilon}$ be one of the connected components of $V_{i,\bar \varepsilon}$. From Lemma \ref{second.index.thm.lemma2} we obtain the existence of $T_i\in\mathcal S_0$ so that
 $\Phi_i(\tilde V_{i,\bar \varepsilon})\subset {\bf B}^{\mathcal F}_{\eta/C}(T_i)$ and thus 
\begin{equation}\label{second.index.thm.flat}
\mathcal F((F_{v})_{\#}\Phi_i(x), (F_{v})_{\#}T_i)< \eta \quad \mbox{for all }x\in \tilde V_{i,\bar \varepsilon}, v\in J.
\end{equation}
In particular this holds for $x \in \partial \tilde V_{i,\bar \varepsilon}$.
Therefore, we can invoke Proposition \ref{f.close}  and derive the existence of  a continuous map $H^2_i: \partial \tilde V_{i,\bar \varepsilon}\times [0,1]\rightarrow \mathcal{Z}_n(M;{\bf M};G)$   such that
\begin{itemize}
\item[(a)] $H^2_i(x,0)=(F_{D_i(x,1)})_{\#}\Phi_i(x)$, $H^2_i(x,1)=(F_{D_i(x,1)})_{\#}T_i$;
\item[(b)]  $H^2_i(x,[0,1])\subset {\bf B}^{\mathcal F}_{\tilde \delta}((F_{D_i(x,1)})_{\#}T_i); $
\item[(c)] ${\bf M}(H^2_i(x,t))\leq {\bf M}((F_{D_i(x,1)})_{\#}\Phi_i(x))+\tilde\delta$ for all $0\leq t\leq 1$.
\end{itemize}
From the third property and Lemma \ref{T.lemma}  we obtain that for all $x\in \partial \tilde V_{i,\bar \varepsilon}$ and $0\leq t\leq 1$
\begin{equation}\label{mass.bound.theorem}
{\bf M}(H^2_i(x,t))\leq  {\bf M}(\Phi_i(x))-\frac{c_0}{10}+\tilde\delta\leq {\bf M}(\Phi_i(x))-3\frac{c_0}{40}.
\end{equation}

If $|T_i|\neq \Sigma$, we have from property ($h_0$) that 
$$ {\bf M}((F_{v})_{\#}T_i)\leq {\bf M}(\Sigma)-c_0\quad\mbox{for all }v\in J.$$ Hence there exists a continuous map $\Xi_i:\tilde V_{i,\bar \varepsilon}\rightarrow \mathcal S_1$ so that $\Xi_i$ equals $H^2_i(\cdot,1)$ on $\partial \tilde V_{i,\bar \varepsilon}$, $\Xi_i(\tilde V_{i,\bar \varepsilon}) \subset \{(F_{v})_{\#}T_i;v\in J\}$ and 
\begin{equation}\label{mass.bound2.theorem}
\sup_{x\in \tilde V_{i,\bar \varepsilon} }{\bf M}(\Xi_i(x))\leq ||\Sigma||(M)-c_0.
\end{equation}

If $|T_i|=|\Sigma|$, then $T_i=\sum_p n_p\cdot S_p$ with $|n_p|=|m_p|$. Because $m_1\geq 2$ we must necessarily have $G=\mathbb{Z}$ and $S_1$ orientable.  By possibly changing the orientation of $S_1$, we can assume $n_1=m_1\geq 2$. We proceed as follows.  Set 
$$A=\{(u,v)\in J\times J: |u|\geq c_0/4 \mbox{ or } |v|\geq c_0/4\} $$
and write $\tilde V_{i,\bar \varepsilon}=[a_i,b_i]$, where $\partial   \tilde V_{i,\bar \varepsilon} = \{a_i,b_i\}$. Denote by $\gamma_i(t)=(u_i(t),v_i(t))$, $a_i\leq t\leq b_i$, a path in $A$ with $\gamma_i(a_i)=(D_i(a_i,1), D_i(a_i,1))$ and $\gamma_i(b_i)=(D_i(b_i,1), D_i(b_i,1))$.

Note that $(F_v)_{\#}(T_i)=m_1 \cdot (F_v)_{\#}(S_1)+\sum_{p\geq 2} n_p\cdot S_p$.  Define $\Xi_i:\tilde V_{i,\bar \varepsilon}\rightarrow \mathcal Z_{n}(M;{\bf F};G)$ as
$$\Xi_i(t)=(F_{u_i(t)})_{\#}S_1+(m_1-1)(F_{v_i(t)})_{\#}S_1+\sum_{p\geq 2} n_p\cdot S_p.$$
The map $\Xi_i$ equals $H^2_i(\cdot,1)$ on $\partial \tilde V_{i,\bar \varepsilon}$ and
\begin{eqnarray}\label{mass.bound3.theorem}
\nonumber{\bf M}(\Xi_i(t))&\leq& {\bf M}((F_{u_i(t)})_{\#}S_1)+(m_1-1) {\bf M}(F_{v_i(t)})_{\#}S_1) + \sum_{p\geq 2} |m_p|\cdot {\bf M}(S_p)\\
\nonumber &\leq& |S_1|-\frac{c_0}{2m_1}|u(t)|^2+(m_1-1)|S_1|-(m_1-1)\frac{c_0}{2m_1}|v(t)|^2\\
\nonumber&&+ \sum_{p\geq 2} |m_p|\cdot {\bf M}(S_p)\\
&\leq& ||\Sigma||(M)-\frac{c_0^3}{32}
\end{eqnarray}
because  $\max\{|u(p,t)|, |v(p,t)|\}\geq c_0/4$. We have used the inequality (\ref{area.inequality}) applied to $V=\Sigma$ and the fact that the diffeomorphisms $F_v$ differ from the identity only in a small neighborhood of $S_1$.

We define  a map  $\bar\Psi_i$ on the space  $\left(\partial \tilde V_{i,\bar \varepsilon}\times[0,1]\right)\cup \left(\tilde V_{i,\bar \varepsilon}\times\{1\}\right)$ by
\begin{eqnarray*}\label{Phimap} 
 \bar \Psi_i(w)=\left\{
\begin{array}{rl}
 (F_{H_i^1(x,3t)})_{\#}\Phi_i(x)& \mbox{if  } w=(x,t)\in  \partial \tilde V_{i,\bar \varepsilon}\times[0,1/3],\\
(F_{D_i(x,3t-1)})_{\#}\Phi_i(x) & \mbox{if  } w=(x,t)\in  \partial \tilde V_{i,\bar \varepsilon}\times[1/3,2/3],\\
H^2_i(x,3t-2) & \mbox{if  } w=(x,t)\in  \partial \tilde V_{i,\bar \varepsilon}\times[2/3,1],\\
 \Xi_i(x) & \mbox{if  } w=(x,1)\in   \tilde V_{i,\bar \varepsilon}\times\{1\}.\\
\end{array}
\right.
\end{eqnarray*}
The map $\Psi_i:\tilde V_{i,\bar \varepsilon}\rightarrow \mathcal Z_{n}(M;{\bf F};G)$ is defined by composing $\bar \Psi_i$ with a  
homemomorphism from $\tilde V_{i,\bar \varepsilon}$ to $\left(\partial \tilde V_{i,\bar \varepsilon}\times[0,1]\right)\cup \left(\tilde V_{i,\bar \varepsilon}\times\{1\}\right)$. 

Since $\Psi_i(y)=\Phi_i(y)$ for every $y \in \partial \tilde V_{i,\bar \varepsilon}$, we can apply the constructions above to every connected component of $V_{i,\bar \varepsilon}$ and get a well-defined continuous map $\Psi_i: X \rightarrow \mathcal Z_{n}(M;{\bf F};G)$ with $\Psi_i=\Phi_i$ on $\overline{X \setminus V_{i, \bar \varepsilon}}$. 

\medskip

\noindent{\bf Claim 1:}  $\Psi_i$ is homotopic to $\Phi_i$  in the flat topology for all $i\in\N$.
\medskip

It suffices to show that $(\Psi_i)_{|\tilde V_{i,\bar \varepsilon}}$ is homotopic to $(\Phi_i)_{|\tilde V_{i,\bar \varepsilon}}$ relative to $\partial \tilde V_{i,\bar \varepsilon}$ in the flat topology for all $i\in\N$ and every connected component $\tilde V_{i,\bar \varepsilon}$.
From  \eqref{second.index.thm.flat} and property (b) of $H^2_i$  we have that 
$$\bar \Psi_i\left((\partial \tilde V_{i,\bar \varepsilon}\times[0,1])\cup (\tilde V_{i,\bar \varepsilon}\times\{1\})\right)\subset {\bf B}^{\mathcal F}_{2\delta}(T_i).$$ Recall that 
$\Phi_i(\tilde V_{i,\bar \varepsilon})\subset {\bf B}^{\mathcal F}_{\eta/C}(T_i)\subset {\bf B}^{\mathcal F}_{\delta/2}(T_i)$ and so the result follows from the choice of $\delta$.
\medskip

\noindent{\bf Claim 2:}  $\limsup_{i\to\infty}\sup\{{\bf M}(\Psi_i(x)):x\in   V_{i,\bar \varepsilon}\}\leq L$.
\medskip

The image of $H^1_i$ is contained in $B_{1/i}^1(0)$ and thus, because $||F_v-\text{Id}||_{C^2}$ tends to zero uniformly as  $v\in J$ tends to zero, we have that
\begin{multline*}\label{mass.inequality.index}
\limsup_{i\to\infty}\sup\{{\bf M}(\bar \Psi_i(y)):y\in \partial V_{i,\bar \varepsilon}\times[0,1/3]\}\\
=\limsup_{i\to\infty}\sup\{{\bf M}((F_{H^1_i(x,t)})_{\#}(\Phi_i(x))):(x,t)\in \partial  V_{i,\bar \varepsilon}\times[0,1]\}\leq L.
\end{multline*}

From the mass non-increasing property of $\phi^x_i$, we have that for all $(x,t)\in \partial  V_{i,\bar \varepsilon}\times[0,1]$,
\begin{equation}\label{mass.decreasing.theorem}
||(F_{D_i(x,t)})_{\#}(\Phi_i(x))||(M)\leq ||(F_{H^1_i(x,1)})_{\#}(\Phi_i(x))||(M)
\end{equation}
and so
$$\limsup_{i\to\infty}\sup\{{\bf M}(\bar \Psi_i(y)):y\in \partial V_{i,\bar \varepsilon}\times[1/3,2/3]\}\leq L.$$
Thus  the claim follows from \eqref{mass.bound.theorem}, \eqref{mass.bound2.theorem}, and \eqref{mass.bound3.theorem}.

\medskip

\noindent{\bf Claim 3: }{\em There exists $0<\xi< \min \{c_0/40,c_0^3/32 \}$ so that $|\Psi_i|(V_{i,\bar \varepsilon})\cap \overline {\bf B}^{\bf F}_{\xi}(\Sigma)=\emptyset$
for all $i$ sufficiently large.}\\

Again because the image of $H^1_i$ is contained in $B_{1/i}^1(0)$,  
we have that
$$\limsup_{i\to\infty}\{{\bf F}(|\Phi_i(x)|,(F_{H^1_i(x,t)})_{\#}(|\Phi_i(x)|)) :(x,t)\in \partial  V_{i,\bar \varepsilon}\times[0,1]\}=0. $$  Hence from (\ref{distance.boundary}) we have that, for all $i$ sufficiently large,
${\bf F}(|\bar \Psi_i(y)|,\Sigma)\geq \bar \varepsilon/2$ for all $y\in \partial  V_{i,\bar \varepsilon}\times[0,1/3]$.

From \eqref{mass.decreasing.theorem} we have that, for all  sufficiently large $i$, 
 $$||(F_{D_i(x,t)})_{\#}(\Phi_i(x))||(M)\leq ||\Phi_i(x)||(M)+\bar\eta\quad\text{for all }(x,t)\in \partial  V_{i,\bar \varepsilon}\times[0,1],$$
where $\bar\eta=\bar \eta(\bar \varepsilon, \Sigma, \{F_v\})>0$ is given by Lemma \ref{away.lemma}. Hence, we obtain from Lemma \ref{away.lemma} and (\ref{distance.boundary}) that ${\bf F}(|\bar \Psi_i(y)|,\Sigma)\geq 2\bar\eta$ for all $y\in \partial  V_{i,\bar \varepsilon}\times[1/3,2/3]$.

We have by (\ref{mass.bound.theorem}) and property $(h_3)$ that $$
{\bf M}(H^2_i(x,t))\leq {\bf M}(\Phi_i(x))-3\frac{c_0}{40} \leq {\bf M}(\Sigma)-c_0/40.
$$
for all $(x,t)\in \partial  V_{i,\bar \varepsilon}\times[0,1]$. Moreover (independently if $|T_i|\neq \Sigma$ or $|T_i|=\Sigma$ in each connected component) we have
$$
\sup_{x\in V_{i,\varepsilon} }{\bf M}(\Xi_i(x))\leq ||\Sigma||(M)-\min \{c_0,c_0^3/32 \}.
$$

Thus,  we have that ${\bf F}(|\Psi_i(y)|,\Sigma)\geq \bar \gamma$ for all $y$ in $\left(\partial  V_{i,\bar \varepsilon}\times[2/3,1]\right)\cup \left( V_{i,\bar \varepsilon}\times{1}\right)$ and some $\bar \gamma =\bar \gamma(\Sigma, c_0)>0$. This concludes the proof of Claim 3.

\medskip

\noindent{\bf Claim 4: }{\em For all sufficiently large $i$, $|\Psi_i|(V_{i,\bar \varepsilon})\cap K\cap\{V:||V||(M)\geq L-\xi\}=\emptyset$.}\\

We see that if $x\in V_{i,\bar \varepsilon}$ and ${\bf M}(\Psi_i(x))\geq L-\xi$, then  $\Psi_i(x)=(F_{v})_{\#}V$ for some $v\in J$, $V\in \overline {\bf B}^{\bf F}_{2\bar \varepsilon}(\Sigma)$. The claim then follows from  $(h_2)$.

\medskip

This finishes the proof of Deformation Theorem B.  

\end{proof}


\subsection{Deformation Theorem C}\label{third.index.theorem}{\em Suppose the sequence $\{\Phi_i\}$ is pulled-tight, i.e., every varifold $V\in {\bf \Lambda}(\{\Phi_i\})$  with $||V||(M)=L$ is stationary. Let $\{\Sigma^{(1)},\dots,\Sigma^{(Q)}\}\subset \mathcal V_n(M)$ be a collection of stationary integral varifolds such that:
\begin{enumerate}
\item[(a)]  the support of $\Sigma^{(q)}$ is a strictly stable two-sided closed embedded minimal hypersurface $S^{(q)}$ for every $1\leq q \leq Q$; 
\item[(b)] $L=||\Sigma^{(q)}||(M)$ for every $1\leq q \leq Q$.
\end{enumerate}
Then there exist $\xi>0$, $j_0\in\N$ so that  for all $i\geq j_0$ one can find   $\Psi_i: X\rightarrow \mathcal Z_n(M;{\bf F};G)$ such that
\begin{itemize}
\item[(i)] $\Psi_i$ is homotopic to $\Phi_i$  in the flat topology;
\item[(ii)] $\limsup_{i\to\infty}\sup\{{\bf M}(\Psi_i(x)):x\in  X\}\leq L$;
\item[(iii)]  $${\bf \Lambda}(\{\Psi_i\})\subset \left({\bf \Lambda}(\{\Phi_i\})\setminus \cup_{q=1}^Q {\bf B}_\xi^{\bf F}(\Sigma^{(q)})\right) \cup \left(\mathcal V_n(M) \setminus {\bf B}_\xi^{\bf F}(\Gamma)\right),$$
where $\Gamma$ is the collection of all stationary integral varifolds  $V\in \mathcal V_n(M)$ with $L=||V||(M)$.
\end{itemize}
}

\begin{proof}
We are going to consider the case $Q=1$. The general case follows similarly with simple modifications.

Let $\Sigma=\Sigma^{(1)}, S=S^{(1)}$. 
Let $S_1,\dots,S_P$ be the connected components of $S$. Then, as an integral varifold, $$\Sigma = m_1\cdot |S_1| + \cdots + m_P \cdot |S_P|$$ with $m_p\in \mathbb{N}$.

Let $\mathcal S_0$ be the finite set of all $T=\sum_p a_p \cdot S_p \in \mathcal Z_n(M;G)$ with $a_p\in G$ and $|a_p|\leq m_p$. If $G=\mathbb{Z}$ and $S_p$ is non-orientable, we set $a_p=0$.
We choose a constant $\delta>0$
so that for all $T\in \mathcal S_0$, the balls ${\bf B}^{\mathcal F}_{2\delta}(T)$ are mutually disjoint and any two paths  in ${\bf B}^{\mathcal F}_{2\delta}(T)$  with common endpoints are homotopic (not
necessarily in ${\bf B}^{\mathcal F}_{2\delta}(T)$) to each other relative to  the boundary in the flat topology.
The existence of such constant follows from the work of Almgren \cite{almgren} (see also Proposition 3.5 or Appendix 2 of \cite{marques-neves-infinitely}). 

The fact that $S$ is strictly stable and two-sided implies that a neighborhood of $S$ can be foliated by hypersurfaces whose mean curvature vectors point towards $S$. This can be done by using  the first eigenfunction of  the Jacobi operator of $S$. As a consequence, we can define a real function $f$ defined on a neighborhood of $S$ with $\nabla f\neq 0$, such that   $S(s)=f^{-1}(s)$ are smooth embedded hypersurfaces for all $|s|\leq 1$, $S_0=S$, and   $S(s)$ is strictly mean convex for all $0<|s|\leq 1$, i.e.,  $s\langle \nabla f,\vec H(S(s))\rangle <0$.  Denoting by $X=\nabla f/|\nabla f|^2$, we consider the  diffeomorphism onto its image
$$\phi:S\times[-1,1]\rightarrow M\quad\mbox{such that } \frac{\partial \phi}{\partial s}(x,s)=X(\phi(x,s)),$$
where $\phi(S,s)=S(s)$ for all $|s|\leq 1$.

Set $\Omega_r=\phi(S\times(-r,r))$ and  consider
$$P:\Omega_1\times[0,1]\rightarrow \Omega_1 \quad\mbox{such that } P(\phi(x,s),t)=\phi(x,(1-t)s).$$
Write $P_t(x)=P(x,t)$.

\subsection{Proposition}\label{projection}{\em There exists $r_0>0$ such that  $P:\Omega_{r_0}\times [0,1]\rightarrow \Omega_{r_0}$ satisfies
\begin{itemize}
\item[(i)] $P(x,0)=x$ for all $x\in \Omega_{r_0}$, $P(x,t)=x$ for all $x\in S$ and $0\leq t\leq 1$;
\item[(ii)] $P(\Omega_{r},t)\subset \Omega_r$ for all $0\leq t\leq 1$, $r\leq r_0$, and $P(\Omega_{r_0},1)=S$;
\item[(iii)] the map $P_t:\Omega_{r_0}\rightarrow \Omega_{r_0}$, $x\mapsto P(x,t)$, is a diffeomorphism onto its image for all $0\leq t<1$;
\item[(iv)] for all $V\in \mathcal V_n(\Omega_{r_0})$ and every connected component $\Omega$ of $\Omega_{r_0}$, the function $t\mapsto ||(P_t)_{\#}V||(\Omega)$ has strictly negative derivative  unless ${\rm supp}\,V\cap \Omega\subset S\cap \Omega$, in which case it is constant.
\end{itemize}
}

\begin{proof}
Directly from the definition we get that for every $r\leq 1$ the map $P$ restricted to $\Omega_{r}\times[0,1]$ satisfies conditions (i), (ii), and (iii) of Proposition \ref{projection}.

For every $0\leq t<1$, consider the vector field  $Z_t=-\frac{f}{1-t}X$.  Given $y=\phi(x,s)\in \Omega_1$, we have $P_t(y)=P(y,t)=P(\phi(x,s),t)=\phi(x,(1-t)s)$. Hence $f(P_t(y))=(1-t)s$ and
\begin{eqnarray*}
\frac{d}{dt}P_t(y)=-s\frac{\partial \phi}{\partial s}(x,(1-t)s)=-\frac{f(P_t(y))}{(1-t)}X(P_t(y))=Z_t(P_t(y)).
\end{eqnarray*}

Given $x\in \Omega_1$ and an $n$-dimensional subspace $\sigma \subset T_x\Omega_1$,  we use the notation
$$\d_\sigma Z_t(x)=\sum_{i=1}^n\langle D_{v_i} Z_t,v_i \rangle\mbox{ where } \{v_i\}_{i=1}^n\mbox{ is an orthonormal basis for }\sigma.$$   The first variation formula says that for all $V\in \mathcal V_n(\Omega_{1})$,
$$\frac{d}{dt}||(P_t)_{\#}V||(M)=\int\d_{\sigma} Z_t \, d\mu_{W}(x,\sigma), \quad\mbox{where }W=(P_t)_{\#}V.$$

Given $y=\phi(x,s)\in\Omega_{1}$,  we can choose local coordinates $(x_1,\ldots,x_n)$ in  a neighborhood of $x\in S$. Then $e_i=\frac{\partial\phi}{\partial x_i}(x,s)$, $i=1,\ldots, n$, span $T_yS(s)$. Let $N=\nabla f/|\nabla f|$, and let $A$ be the second fundamental form of $S(s)$. We have 
\begin{align*}
\langle D_{e_i} Z_t,e_j\rangle & =-\langle Z_t, A(e_i,e_j)\rangle,\\
\langle D_{e_i} Z_t,N\rangle & =-\frac{f}{1-t}\left\langle D_{\frac{\partial\phi}{\partial x_i}}\frac{\partial\phi}{\partial s},
 N\right\rangle=\frac{f}{1-t}\left\langle D_{\frac{\partial\phi}{\partial s}} N, \frac{\partial\phi}{\partial x_i}\right \rangle
= -\langle D_{N} Z_t,e_i\rangle,\\
 \langle D_{N} Z_t,N\rangle & =-\frac{1}{1-t}(1+f\langle D_{N} X,N\rangle).
\end{align*}

Choose   an orthonormal basis $\{v_i\}_{i=1}^n$ for $\sigma$ so that $v_1,\ldots, v_{n-1}$ are orthogonal to $N$. 
Then $v_n=\cos\theta u+\sin\theta N$ for some $u\in T_yS(s)$, $|u|=1$, and angle $\theta$. The expressions above show that $\sum_{i=1}^n\langle D_{v_i} Z_t,v_i\rangle$ is equal to
$$
\frac{1}{1-t}\left(|\nabla f|^{-1} s\langle N,\vec H(S_l)\rangle-\sin^2\theta(1+s\langle D_{N} X,N\rangle+s\langle X, A(u,u)\rangle)\right).
$$
Thus there is $r_0$  so that if $|s|\leq r_0$ we obtain that $\sum_{i=1}^n\langle D_{v_i} Z_t,v_i\rangle$  is always non-positive and zero only if $s=0$. This shows property (iv).

\end{proof}

We deduce the following corollary:

\subsection{Corollary}\label{projection.corollary}{\em There exists $\varepsilon_0=\varepsilon_0(\Sigma)>0$ so that every stationary integral varifold $V\in  \mathcal V_n(M)$ in ${\bf B}^{\bf F}_{\varepsilon_0}(\Sigma)$ coincides with $\Sigma$.
}

\medskip

The above result does not hold if one does not require the stationary varifold $V$ to  be integral.
\begin{proof}
Suppose there is a sequence of   stationary integral varifolds $\{V_j\}_{j\in\N}$ in $\mathcal V_n(M)$ converging to $\Sigma$ in the varifold norm with $V_j \neq \Sigma$. Thus $||V_j||(M\setminus \Omega_{r_0/2})$ tends to zero as $j$ tends to infinity and so the monotonicity formula implies that ${\rm supp}\, V_j\subset \Omega_{r_0}$ for all $j$ sufficiently large. From Proposition \ref{projection} (iv) we deduce that ${\rm supp}\, V_j\subset S$, which when combined with the Constancy Theorem for stationary varifolds implies that $V_j=\sum_p n_{j,p} \cdot |S_p|$ for $n_{j,p} \in \mathbb{Z}_+$. Since ${\bf F}(V_j,\Sigma) \rightarrow 0$, we necessarily have $V_j=\Sigma$ for all $j$ sufficiently large. Contradiction.
\end{proof}

From Proposition \ref{F.approximation}  we  see that we can assume, without loss of generality, that  $\Phi_i$ 
 is continuous in the mass topology for every $i$.

Let $X(k_i)$ be a sufficiently fine subdivision of $X$ so that ${\bf M}(\Phi_i(x)-\Phi_i(y))<1/i$ and ${\bf F}(|\Phi_i(x)|,|\Phi_i(y)|)<1/i$ for any $x,y$ belonging to the
 same 1-cell in $X(k_i)$. For the purpose of this proof we redefine $V_{i,\alpha}$ to be the union of all 1-cells $\sigma \in X(k_i)$ so that ${\bf F}(|\Phi_i(x)|,\Sigma)<2\alpha$ for all $x\in \sigma$. 
 If $\beta \notin V_{i,\alpha}$, then ${\bf F}(|\Phi_i(x')|,\Sigma)\geq 2\alpha$ for some $x'\in \beta$. Consequently,  ${\bf F}(|\Phi_i(x)|,\Sigma)\geq 2\alpha-1/i$ for all $x\in \beta$.  

\subsection{First construction:}

The next Proposition deforms those cycles $\Phi_i(x)$ that are very close to $\Sigma$ in varifold norm to cycles supported in a tubular neighborhood of $S$. We use the deformation given by Lemma 7.1 of Montezuma \cite{montezuma}.

\subsection{Proposition}\label{first.construction}{\it There exists  a constant $C_1=C_1(M)>0$ so that for all $\eta>0$ one can find  $\tilde \varepsilon=\tilde \varepsilon(\Sigma, r_0, \eta)>0$  with the property that for all $0<\varepsilon\leq \tilde \varepsilon$, there exists  $i_1\in\N$ such that if $i\geq i_1$, there is  a continuous map  in the flat topology
$$H^1_i: V_{i,\varepsilon}\times[0,1] \rightarrow \mathcal Z_n(M;\mathcal F;G)$$
such that
\begin{itemize}
\item[(i)] $H^1_i(x,0)=\Phi_i(x)$ and ${\rm supp}\,H^1_i(x,1)\subset \Omega_{r_0/2}$ for all $x\in  V_{i,\varepsilon}$;
\item[(ii)] for all $(x,t)\in  V_{i,\varepsilon}\times[0,1]$  we have  ${\bf M}(H^1_i(x,t)-\Phi_i(x))\leq C_1\eta $ and in particular ${\bf M}(H^1_i(x,t))\leq {\bf M}(\Phi_i(x))+C_1\eta$;
\item[(iii)] for all $(x,t)\in \partial  V_{i,\varepsilon}\times[0,1]$ we have ${\bf M}(H^1_i(x,t))\leq {\bf M}(\Phi_i(x))+C_1/i$ and $H^1_i(x,t)\llcorner \Omega_{r_0/4}=\Phi_i(x)\llcorner \Omega_{r_0/4}$;
\item[(iv)] $H^1_i$ restricted to $(\partial V_{i,\varepsilon} \times [0,1]) \cup (V_{i,\varepsilon}\times \{1\})$ is continuous in the mass topology.
\end{itemize}}

\begin{proof}
Let  $U=M\setminus \bar \Omega_{r_0/3}$ and $W=M\setminus \bar \Omega_{2r_0/5}$. Consider  
$\varepsilon_1=\varepsilon_1(U,W)>0$ given by \cite[Lemma 7.1]{montezuma}.
In what follows let $\eta>0$ be such that $5\eta<\min\{\delta_0,\varepsilon_1\}$, where $\delta_0$ is the  constant provided by Theorem 3.10 of \cite{marques-neves-infinitely} (the same statement holds if $G=\mathbb{Z}$).
 There exists  $\bar \varepsilon=\bar \varepsilon(\Sigma,r_0,\eta)>0$ so that 
$${\bf F}(|\Phi_i(x)|,\Sigma)\leq  2\bar \varepsilon \implies ||\Phi_i(x)||(U)< \eta.$$ 
We  assume that $0< \varepsilon\leq \min\{\bar \varepsilon, \eta\}.$

Applying  Lemma 7.1 of \cite{montezuma} to each image of $\Phi_i$ restricted to $(V_{i,\varepsilon})_0$  and observing that the deformations are local, we obtain $\hat H^1_i:(V_{i,\varepsilon})_0\times I(1,l_i)_0\rightarrow \mathcal  Z_n(M;G)$ for some $l_i\in\N$ sufficiently large such that
\begin{itemize}
\item[(a1)] $\hat H^1_i(x,0)=\Phi_i(x)$ and ${\rm supp}\,\hat H^1_i(x,1)\subset \Omega_{2r_0/5}$ for all $x\in (V_{i,\varepsilon})_0$;
\item[(b1)] ${\rm supp} (\hat H_i^1(x,t)-\Phi_i(x)) \subset U$ for all $(x,t)\in (V_{i,\varepsilon})_0\times I(1,l_i)_0$;
\item[(c1)] ${\bf M}({\hat H_i^1}(x,t)-{\hat H_i^1}(x,t'))<1/i$ if $(x,t), (x,t')$ are adjacent vertices in  $(V_{i,\varepsilon})_0\times I(1,l_i)_0$;
\item[(d1)] ${\bf M}(\hat H^1_i(x,t)\llcorner U)\leq {\bf M}(\Phi_i(x)\llcorner U)+1/i$ for all $(x,t)\in (V_{i,\varepsilon})_0\times I(1,l_i)_0$.
  \end{itemize}

We have 
\begin{eqnarray*}
{\bf M}(\hat H^1_i(x,t)-\Phi_i(x))&\leq& {\bf M}(\hat H^1_i(x,t)\llcorner U)+{\bf M}(\Phi_i(x)\llcorner U)\\
&\leq& 2\,{\bf M}(\Phi_i(x)\llcorner U)+1/i\\
&\leq& 2\eta+1/i,
\end{eqnarray*}
 for all $(x,t)\in (V_{i,\varepsilon})_0\times I(1,l_i)_0$. Hence,  for all $(x,t), (x',t)$  adjacent vertices in  $(V_{i,\varepsilon})_0\times I(1,l_i)_0$ and $i$ sufficiently large, we have
 \begin{eqnarray*}
 {\bf M}(\hat H^1_i(x,t)-\hat H^1_i(x',t))&\leq& {\bf M}(\hat H^1_i(x,t)-\Phi_i(x)) \\
 &&+ {\bf M}(\Phi_i(x)-\Phi_i(x'))+ {\bf M}(\hat H^1_i(x',t)-\Phi_i(x'))\\
 &\leq&4\eta +3/i\\
 &<& 5\eta.
 \end{eqnarray*}

This implies the  map $\hat H^1_i$ has fineness smaller than $5 \eta<\delta_0$. 
Using Theorem 3.10 of \cite{marques-neves-infinitely} (the same statement holds for $G=\mathbb{Z}$), we obtain a continuous map $\bar H^1_i:V_{i,\varepsilon}\times[0,1] \rightarrow \mathcal Z_n(M;{\bf M};G)$  in the mass topology (the {\it Almgren extension} of $\hat H^1_i$) and $C_0=C_0(M)>0$ such that
\begin{itemize}
\item[(a2)] $\bar H^1_i(x,t)=\hat H^1_i(x,t)$ if $(x,t) \in (V_{i,\varepsilon})_0\times I(1,l_i)_0$ and ${\rm supp}\,\bar H^1_i(x,1)\subset \Omega_{r_0/2}$ for all $x\in V_{i,\varepsilon}$;
\item[(b2)] 
$\lim_{i\to\infty}\sup\{{\bf M}(\Phi_i(x)-\bar H^1_i(x,0)):x\in V_{i,\varepsilon}\}=0;$
\item[(c2)]  ${\bf M}(\bar H^1_i(x,t)-\Phi_i(x))\leq  (5C_0+3)\eta$ for all $(x,t)\in V_{i,\varepsilon}\times I$;
\item[(d2)]  ${\bf M}(\bar H^1_i(x,t))\leq {\bf M}(\Phi_i(x))+(C_0+1)/i$  and $$\bar H^1_i(x,t)\llcorner \Omega_{r_0/4}=\Phi_i(x)\llcorner\Omega_{r_0/4}$$ for all $(x,t)\in \partial V_{i,\varepsilon}\times I$.
\end{itemize}
Almgren interpolation results use a triangulation of $M$. The claims about the supports above follow by choosing the triangulation sufficiently fine with respect to $r_0$ (see proof of Theorem 14.1 of \cite{marques-neves}).

Note that $\bar H^1_i(x,0)=\Phi_i(x)$ for $x\in \partial V_{i,\varepsilon}$. Using Proposition \ref{mass.close} of the Appendix and property (b2),  we find a homotopy of ${\Phi_i}_{|V_{i,\varepsilon}}$ and $\bar H^1_i(\cdot,0)$. We can glue that homotopy with $\bar H^1_i$ to construct  $H^1_i$ satisfying properties (i), (ii), (iii) and (iv).

\end{proof}

\subsection{Second construction:}  
 From the Constancy Theorem, there is $\eta_1=\eta_1(\Sigma)>0$ so that  any $T\in \mathcal Z_n(M;G)$ with ${\rm supp}\, T\subset S_p$ for some $1\leq p\leq P$ and ${\bf M}(T)\leq m_p ||S_p||(M)+2\eta_1$ must satisfy $T=m\cdot S_p$ with $m \in G$ and $|m|\leq m_p$.
 
From the Constancy Theorem for stationary varifolds, there is also $\alpha>0$ (depending only on $\Sigma$) so that for every stationary integral varifold $Z$ with ${\rm supp}\, Z \subset S$, ${\bf M}(Z)={\bf M}(\Sigma)$, but $\Sigma\neq Z$, there is $\Omega^Z$ a connected component of $\Omega_{r_0}$ so that $||Z||(\Omega^Z)>||\Sigma||(\Omega^Z)+2\alpha.$ Moreover there is $\beta>0$ (depending on $\Sigma$ and $\alpha$) so that every $W \in {\bf B}^{\bf F}_{\beta}(\Sigma)$ has $||W||(\Omega)\leq ||\Sigma||(\Omega)+\alpha/2$ for every connected component $\Omega$ of $\Omega_{r_0}$.

Consider $\varepsilon_0>0$ given by Corollary \ref{projection.corollary}.   We choose $\eta>0$ (as in Proposition \ref{first.construction}) and $\varepsilon>0$  so that $\eta <\min\{(2\eta_1)/(C_1+1), \varepsilon_0/(6C_1)\}$, $\varepsilon<\min\{\tilde{\varepsilon}(\Sigma, r_0,\eta),\varepsilon_2/2,\varepsilon_0/4, \beta/2\}$, where $\varepsilon_2=\varepsilon_2(\delta,\Sigma)$ is as in 
Lemma \ref{second.index.thm.lemma2}. We can also require
$$
{\bf F}(V,\Sigma)\leq 2\varepsilon \implies ||V||(\Omega_{r_0/2}^{(p)})\leq m_p ||S_p||(M)+ \eta,
$$
where $\Omega_{r_0/2}^{(p)}$ denotes the connected component of $\Omega_{r_0/2}$ that contains $S_p$, $1\leq p\leq P$.

Like in the proof of Deformation Theorem B we can assume that $|\Phi_i|(X)$ is not contained in ${\bf B}^{\bf F}_{\varepsilon_2}(\Sigma)$ and that intersects ${\bf B}^{\bf F}_{\varepsilon}(\Sigma)$. Thus  $V_{i,\varepsilon}$ has finitely many connected components which are homeomorphic to closed intervals and  $\partial V_{i,\varepsilon}$ is a non-empty union of vertices of $X(k_i)$.

 For all $i\geq i_1$ consider  the continuous map in the flat topology
$$H^2_i:V_{i,\varepsilon}\times[0,1] \rightarrow \mathcal Z_n(M;G),\quad H^2_i(x,t)=(P_t)_{\#}(H^1_i(x,1)).$$

Denote by $\tilde V_{i,\varepsilon}$ a connected component of $V_{i,\varepsilon}$.
\subsection{Lemma}\label{second.construction}{\em There exists $i_2\in \N$ so that for all $i\geq i_2$:
\begin{itemize}
\item[(i)]${\rm supp}\, H^2_i(x,t)\subset \Omega_{r_0/2}$ for all $(x,t)\in  \tilde V_{i,\varepsilon}\times[0,1]$ and  $H^2_i$ is continuous in the ${\bf F}$-topology  on  $\tilde V_{i,\varepsilon}\times[0,1)$;
\item[(ii)] $ ||H^2_i(x,t)||(\Omega)\leq ||H^1_i(x,1)||(\Omega)$ for every $(x,t)\in  \tilde V_{i,\varepsilon}\times[0,1]$ and every connected component $\Omega$ of $\Omega_{r_0}$;
\item[(iii)]  there exists $\Sigma_i\in \mathcal S_0$ so that $H_i^2(x,1)=\Sigma_i$ for all $x\in \tilde V_{i,\varepsilon}$.
\end{itemize}}

\begin{proof}
The first  and second properties follow from Proposition \ref{projection} and Proposition  \ref{first.construction} (i).

From Proposition \ref{first.construction} (ii) we obtain the existence of $i_2\in\N$ so that for all $i\geq  i_2$ and $1\leq p\leq P$ we have
\begin{eqnarray*}
{\bf M}(H_i^2(x,1) \llcorner \Omega_{r_0/2}^{(p)}) &\leq& {\bf M}(H_i^1(x,1) \llcorner \Omega_{r_0/2}^{(p)})\\
&\leq&{\bf M}(\Phi_i(x)\llcorner \Omega_{r_0/2}^{(p)})+ C_1\eta\\
&\leq&m_p ||S_p||(M)+ (C_1+1)\eta\\
&\leq&m_p ||S_p||(M)+ 2\eta_1,
\end{eqnarray*}
From the choice of $\eta_1$ we conclude that $H_i^2(x,1)\in \mathcal S_0$ for every $x\in \tilde V_{i,\varepsilon}.$
The third property follows from the facts that $\mathcal S_0$ is finite, $H_i^2$ is continuous in the flat topology and $\tilde V_{i,\varepsilon}$ is connected.

\end{proof}

\subsection{Third construction:}\label{third.construction} 
For  $i\in \mathbb{N}$  sufficiently large we have  $1/i<\min\{\delta/P,\nu_3\}$, where $P$ is the number of connected components of $S$ and $\nu_3$ is the constant of Corollary 1.14 of \cite{almgren}.   If
$\mathcal F(T')<\nu_3$ for some $T'\in \mathcal Z_n(M,G)$, then there exists $D'\in I_{n+1}(M,G)$ with $\partial D'=T'$ and ${\bf M}(D')=\mathcal F(T')$. Consider $\eta_i=\eta_i(\mathcal S_0,1/i, L,W,U)$ given by Proposition \ref{f.close}  of the Appendix, where $W=\Omega_{r_0/2}$ and $U=\Omega_{r_0}$.  

Let $V_{i,\varepsilon}^{(l)}$, $l=1,\dots, N_i$, be the connected components of $V_{i,\varepsilon}$ and $\Sigma_i^{(l)}$ be the corresponding cycles given by  Lemma \ref{second.construction} (iii). We obtain the existence of $0<q_i<1/i$ so that $$H^2_i(V_{i,\varepsilon}^{(l)}\times[1-q_i,1])\subset {\bf B}^{\mathcal F}_{\eta_i}(\Sigma_i^{(l)})\quad\mbox{for all }l\leq N_i, i\geq i_2.$$ Note that ${\rm supp}\, H^2_i(x,1-q_i)\subset \Omega_{r_0/2}$ for all $(x,t)\in   V_{i,\varepsilon}\times[0,1]$, and ${\rm supp}\, T \subset S \subset \Omega_{r_0/2}$ for all $T\in \mathcal S_0$.
Then (by the choice of $\nu_3$ and the Constancy Theorem) we have $\mathcal F(H^2_i(x,1-q_i)\llcorner \Omega, \Sigma_i^{(l)}\llcorner \Omega) < \eta_i$ for every connected component $\Omega$ of $\Omega_{r_0}$ and $x\in \partial V_{i,\varepsilon}^{(l)}$.
 Therefore, we can invoke Proposition \ref{f.close}  (since $\Sigma_i^{(l)}\llcorner \Omega \in \mathcal S_0$) and, for all $i\geq i_2$, conclude the existence of  a continuous map $H^3_i: \partial V_{i,\varepsilon}^{(l)}\times [0,1]\rightarrow \mathcal{Z}_n(M;{\bf M};G)$   such that
\begin{itemize}
\item[(a)] $H^3_i(x,0)=H^2_i(x,1-q_i)$, $H^3_i(x,1)=\Sigma_i^{(l)}$;
\item[(b)]  $H^3_i(x,t)\subset {\bf B}^{\mathcal F}_{P/i}(\Sigma_i^{(l)})$ and ${\rm supp}\, H^3_i(x,t)\subset \Omega_{r_0}$ for all $(x,t)\in \partial V_{i,\varepsilon}^{(l)}\times [0,1]$; 
\item[(c)] $||H^3_i(x,t)||(\Omega)\leq ||H^2_i(x,1-q_i)||(\Omega)+1/i$ for every $(x,t)\in \partial V_{i,\varepsilon}^{(l)}\times [0,1]$ and every connected component $\Omega$ of $\Omega_{r_0}.$
\end{itemize}

We set $H^3_i(x,1)= \Sigma_i^{(l)}$ for $x \in V_{i,\varepsilon}^{(l)}$. By the choice of $\delta$, the maps $H_i^2(\cdot,1-q_i)_{|V_{i,\varepsilon}^{(l)}}$ and $(H_i^3)_{|(\partial V_{i,\varepsilon}^{(l)}\times [0,1]) \cup (V_{i,\varepsilon}^{(l)} \times \{1\})}$ are homotopic to each other in the flat topology.

\subsection{Final construction:}

For each $l=1,\dots, N_i$, we define a continuous map in the ${\bf F}$-metric $\bar\Psi_i$ on  $\left(\partial V_{i,\varepsilon}\times[0,1]\right)\cup \left(V_{i,\varepsilon}\times\{1\}\right)$ as
\begin{eqnarray*}\label{Phimap} 
 \bar\Psi_i(w)=\left\{
\begin{array}{rl}
 H_i^1(x,3t) & \mbox{if  } w=(x,t)\in  \partial V_{i,\varepsilon}^{(l)}\times[0,1/3],\\
H_i^2(x,3(1-q_i)t-1+q_i) & \mbox{if  } w=(x,t)\in  \partial V_{i,\varepsilon}^{(l)}\times[1/3,2/3],\\
H^3_i(x,3t-2) & \mbox{if  } w=(x,t)\in  \partial V_{i,\varepsilon}^{(l)}\times[2/3,1],\\
H_i^3(x,1) & \mbox{if  } w=(x,1)\in   V_{i,\varepsilon}^{(l)}\times\{1\}\\
\end{array}
\right.
\end{eqnarray*}
The map $\Psi_i:V_{i,\varepsilon}^{(l)}\rightarrow \mathcal Z_{n}(M;{\bf F};G)$ is defined by composing $\bar \Psi_i$ with a  homeomorphism from $V_{i,\varepsilon}^{(l)}$ to $\left(\partial V_{i,\varepsilon^{(l)}}\times[0,1]\right)\cup \left(V_{i,\varepsilon}^{(l)}\times\{1\}\right)$. Since $\Psi_i(x)=\Phi_i(x)$ for $x\in \partial V_{i,\varepsilon}^{(l)}$, this defines a map $\Psi_i:X \rightarrow \mathcal Z_{n}(M;{\bf F};G)$.

Property (i) of Deformation Theorem C follows from the construction of $\bar \Psi_i$.  Property (ii) follows from  Proposition \ref{first.construction} (iii), Proposition \ref{second.construction} (ii), and property (c) in Subsection \ref{third.construction}.

We now prove Deformation Theorem C (iii).
\medskip


 \noindent{\bf Claim 1:} There exist $\varepsilon_3=\varepsilon_3(\Sigma,r_0,\eta,\{H^1_i\}_i)>0$ and $i_3\in \N$ so that for all $i\geq i_3$ we have $$|H^1_i|(\partial V_{i,\varepsilon}\times[0,1])\cap {\bf B}^{\bf F}_{\varepsilon_3}(V)=\emptyset$$ 
for every stationary integral  varifold $V\in \mathcal V_n(M)$ with $L=||V||(M).$
\medskip

It suffices to show that there exists no subsequence $\{j\}\subset \{i\}$ and  $\{(x_j,t_j)\}_{j}$ in $\partial V_{j,\varepsilon}\times[0,1]$ such that for some stationary integral varifold $V_j$ with $L=||V_j||(M)$ we have  ${\bf F}(|H^1_j(x_j,t_j)|, V_j)\rightarrow 0$ as $j\rightarrow \infty$.  We can further assume $V_j$ converges to $V$ in the ${\bf F}$-metric. 

If such a sequence exists,  we have from  Proposition \ref{first.construction} (ii) and the fact that ${\bf F}(|S|,|T|)\leq {\bf M}(S-T)$ for any $S,T\in \mathcal Z_n(M;G)$ (see Pitts \cite{pitts}, Section 2.1),  $${\bf F}(|H^1_j(x_j,t_j)|,\Sigma)<2\varepsilon+C_1\eta\leq 2\varepsilon_0/3.$$ Hence 
$V_j\in {\bf B}^{\bf F}_{\varepsilon_0}(\Sigma)$ for sufficiently large $j$ and so $V=V_j=\Sigma$ by Corollary \ref{projection.corollary}. Thus
$$\lim_{j\to\infty}{\bf M}(H^1_j(x_j,t_j)\llcorner \Omega_{r_0/8})=L, \quad    \lim_{j\to\infty}{\bf M}(H^1_j(x_j,t_j)\llcorner (M\setminus\Omega_{r_0/8}))=0.$$
We obtain from Proposition \ref{first.construction} (iii) that
$$\lim_{j\to\infty}  {\bf M}(\Phi_j(x_j)\llcorner \Omega_{r_0/8})=L.$$
Since $\limsup_{j\rightarrow \infty} {\bf M}(\Phi_j(x_j))\leq L$, we also get
$$\lim_{j\to\infty}  {\bf M}(\Phi_j(x_j)\llcorner (M\setminus\Omega_{r_0/8}))=0.$$
Therefore
\begin{multline*} \lim_{j\to\infty}  {\bf M}(H^1_j(x_j,t_j)-\Phi_j(x_j)) \leq \lim_{j\to\infty}\left({\bf M}((H^1_j(x_j,t_j)-\Phi_j(x_j))\llcorner \Omega_{r_0/4}) \right.\\
\left.+  {\bf M}(\Phi_j(x_j)\llcorner (M\setminus\Omega_{r_0/8}))+ {\bf M}(H^1_j(x_j,t_j)\llcorner (M\setminus\Omega_{r_0/8}))\right)=0
\end{multline*}
and so $|\Phi_j(x_j)|$ converges to $\Sigma$ in varifold sense as $j$ tends to infinity. This is impossible because $x_j\in \partial  V_{j,\varepsilon}$ for all $j$. 
\medskip

\noindent{\bf Claim 2:} There is $\varepsilon_4=\varepsilon_4(\Sigma,r_0,\eta,\{H^2_i\}_i)>0$ and $i_4\in \N$ so that for all $i\geq i_4$ we have $$|H^2_i|(\partial V_{i,\varepsilon}\times[0,1))\cap {\bf B}^{\bf F}_{\varepsilon_4}(V)=\emptyset$$
for every stationary integral  varifold  $V\in \mathcal V_n(M)$ with $L=||V||(M).$ 
\medskip

Like before, it suffices to show that there exists no subsequence $\{j\}\subset \{i\}$ and  $\{(x_j,t_j)\}_{j}$ in $\partial V_{j,\varepsilon}\times[0,1)$ 
 such that for some stationary integral varifold $V_j$ with $L=||V_j||(M)$ we have  ${\bf F}(|H^2_j(x_j,t_j)|, V_j)\rightarrow 0$ as $j\rightarrow \infty$.  We can further assume $V_j$ converges to $V$ in the ${\bf F}$-metric.

If such a sequence exists, we have from Lemma \ref{second.construction} (i)  that ${\rm supp}\, V\subset \bar \Omega_{r_0/2}$  and thus ${\rm supp}\, V\subset S$ by Proposition \ref{projection} (iv).   
After passing to a subsequence if necessary we can assume  that, as $j$ tends to infinity, $t_j$ tends to $\bar t\in [0,1]$ and $|H^1_j(x_j,1)|$ tends to $\tilde V$, where $||\tilde V||(M)\leq L$. From the definition of $H^2_i$ and  Lemma \ref{second.construction} (ii) we see that  $V=(P_{\bar t})_{\#}(\tilde V)$ and $||V||(M)\leq ||\tilde V||(M)$. Thus we have from Proposition \ref{projection} (i) and (iv)  that $\tilde V=V$. This contradicts Claim 1.

\medskip

\noindent{\bf Claim 3:} There is $\delta_1=\delta_1(\Sigma,r_0,\eta,\{H^2_i\}_i)>0$ and $i_5\in \N$ so that for all $i\geq i_5$, $x\in \partial V_{i,\varepsilon}$  one can find $\Omega_i$, a  connected component of $\Omega_{r_0}$, so that   $$||(P_1)_{\#}(|H^1_i(x,1)|)||(\Omega_{i})\leq ||\Sigma||(\Omega_{i})- 2\delta_1.$$

We argue by contradiction and assume there is a subsequence $\{j\}\subset \{i\}$ and $\{x_j\}_{j}$ in  $ \partial V_{j,\varepsilon}$ so that for every connected component $\Omega$ of $\Omega_{r_0}$ we have
$$\limsup_{j\to\infty}||(P_1)_{\#}(|H^1_j(x_j,1)|)||(\Omega)\geq ||\Sigma||(\Omega).$$
We can assume $(P_1)_{\#}(|H^1_j(x_j,1)|)$  converges, in varifold sense, to $V\in \mathcal V_n(M)$ as $j$ tends to infinity, where $||V||(M)=L$ and ${\rm supp}\, V\subset S$. Arguing like in the proof of Claim 2,  we conclude that $\{|H^1_j(x_j,1)|\}_{j}$  converges to $V\in \mathcal V_n(M)$ as well. Thus
$$\lim_{j\to\infty}{\bf M}(H^1_j(x_j,1)\llcorner (M\setminus\Omega_{r_0/8}))=0\quad\mbox{and}\quad \lim_{j\to\infty}{\bf M}(H^1_j(x_j,1)\llcorner \Omega_{r_0/8})=L.$$
We can argue like in the proof of Claim 1 to conclude that 
$$ \lim_{j\to\infty}  {\bf M}(H^1_j(x_j,1)-\Phi_j(x_j))=0,$$
from which it follows that $\{|\Phi_j(x_j)|\}_{j}$  converges to $V$. Hence ${\bf F}(V,\Sigma)\geq \varepsilon$. 

On the other hand, $||V||(M)=L$ and hence $V$ is stationary because the sequence  $\{{\Phi_i}\}$ is pulled-tight. Since we have 
$$||V||(\Omega)\geq ||\Sigma||(\Omega)\quad\mbox{for every connected component }\Omega\mbox{ of }\Omega_{r_0},$$
the Constancy Theorem for stationary varifolds implies that $V=\Sigma$, which is a contradiction.
\vskip 0.05in

From Claim 3 we see that we can choose $0<q_i<1/i$ in Section \ref{third.construction} so that for every $i\geq i_5$ and  $x\in  \partial V_{i,\varepsilon}$ we have $||H^2_i(x,1-q_i)||(\Omega_i(x))\leq ||\Sigma||(\Omega_i(x))-\delta_1$. 

\medskip

\noindent{\bf Claim 4:} There exist $\varepsilon_5=\varepsilon_5(\Sigma,r_0,\eta,\{H^3_i\}_i)>0$ and $i_6\in \N$ so that for all $i\geq i_6$ we have $$|H^3_i|(\partial V_{i,\varepsilon}\times[0,1])\cap {\bf B}^{\bf F}_{\varepsilon_5}(V)=\emptyset$$
for every stationary integral  $V\in \mathcal V_n(M)$ with $L=||V||(M).$ 
 \medskip

 Like before, suppose there is a subsequence $\{j\}\subset \{i\}$ and $\{(x_j,t_j)\}_{j}$ in $\partial V_{j,\varepsilon}\times[0,1]$  such that for some stationary integral varifold $V_j$ with $L=||V_j||(M)$ we have  ${\bf F}(|H^3_j(x_j,t_j)|, V_j)\rightarrow 0$ as $j\rightarrow \infty$.  We can further assume $V_j$ converges to $V$ in the ${\bf F}$-metric.  Then ${\rm supp}\, V\subset S$. By Allard's compactness theorem \cite{allard} $V$ is also integral.

 From property (c) in Section \ref{third.construction} we have that $$||H^3_j(x_j,t_j)||(\Omega_j(x_j))\leq ||\Sigma||(\Omega_j(x_j))-\delta_1/2$$ for all $j$ sufficiently large. Hence $||V||(\Omega)<||\Sigma||(\Omega)$ for some connected component $\Omega$ of $\Omega_{r_0}$. 
 
 In particular $V\neq \Sigma$ and so, for some connected component $\Omega^V$ of $\Omega_{r_0}$, we have that $||V||(\Omega^V)>||\Sigma||(\Omega^V)+2\alpha.$ Thus, we obtain from Proposition \ref{second.construction} (ii) and property (c) in Section \ref{third.construction} that, for all $j$ sufficiently large, $||H^1_j(x_j,1)||(\Omega^V)>||\Sigma||(\Omega^V)+3\alpha/2.$ By passing to a subsequence, we can assume
 $|H^1_j(x_j,1)|$ converges to $\bar V\in \mathcal V_n(M)$ with $||\bar V||(M) \leq L$, ${\rm supp}\, \bar V \subset \bar \Omega_{r_0/2}$. Since 
 $\lim_{j\rightarrow \infty} ||H^3_j(x_j,t_j)||(M) = L$, we necessarily  have $||\bar V||(M) = L$. It follows that
  $\lim_{j\rightarrow \infty} ||H^2_j(x_j,1/2)||(M) = L$. Since $|H^2_j(x_j,1/2)| \rightarrow (P_{1/2})_{\#}(\bar V)$ as $j \rightarrow \infty$, we have
  $||(P_{1/2})_{\#}(\bar V)||(M)=||\bar V||(M)=L.$ Therefore ${\rm supp}\, \bar V \subset S$. We again have
  $$\lim_{j\to\infty}{\bf M}(H^1_j(x_j,1)\llcorner (M\setminus\Omega_{r_0/8}))=0\quad\mbox{and}\quad \lim_{j\to\infty}{\bf M}(H^1_j(x_j,1)\llcorner \Omega_{r_0/8})=L.$$ Arguing like in Claim 1 we deduce that  ${\bf M}(H^1_j(x_j,1)-\Phi_j(x_j))\rightarrow 0$, so $\bar V$ is stationary because the sequence $\{\Phi_i\}_i$ is pulled-tight. But this also implies that $||\Phi_j(x_j)||(\Omega^V)>||\Sigma||(\Omega^V)+\alpha/2$ for all $j$ sufficiently large. This is impossible because, since $x_j \in V_{j,\varepsilon}$ we have ${\bf F}(\Phi_j(x_j),\Sigma)<2\varepsilon<\beta$.

\smallskip
 
Deformation Theorem C (iii) follows from Claim 1, Claim 2 and Claim 4, together with the fact that $\Psi_i=\Phi_i$ in $X\setminus V_{i,\varepsilon}$ and hence ${\bf F}(\Psi_i(x),\Sigma)\geq \varepsilon$ for all $x\notin V_{i,\varepsilon}$.

\end{proof}

\section{Index bounds}\label{index.bounds}

We assume in this section that $3\leq (n+1)\leq 7$.
Let $X$ be a $k$-dimensional cubical complex and  $\Phi:X \rightarrow \mathcal Z_n(M^{n+1};{\bf F};G)$ be a continuous map. We let $\Pi$ be the class of all continuous maps $\Phi':X \rightarrow \mathcal Z_n(M^{n+1};{\bf F};G)$ such that $\Phi$ and $\Phi'$ are homotopic to each other in the flat topology. 


The metric $(M,g)$ is said to be {\it bumpy} if no smooth immersed closed minimal hypersurface has a non-trivial Jacobi vector field. White showed \cite{white2, white3}   that bumpy metrics are generic in the Baire sense.

\subsection{Theorem}\label{index.k.theorem}{\em 

Assume $(M,g)$ is a bumpy metric and let $\{\Phi_i\}_{i\in\N}$ be a  sequence in $\Pi$ such that $L={\bf L}(\{\Phi_i\}_{i\in\N})={\bf L}(\Pi).$ 

There is $\Sigma\in  {\bf C}(\{\Phi_i\}_{i\in\N})$ with support a smooth embedded closed minimal hypersurface 
 such that
$${\bf L}(\Pi)=||\Sigma||(M)\text{ and index\,(support of $\Sigma$)}\leq k.$$}

\begin{proof}
It suffices to show that, for every $r>0$, there is  a stationary integral varifold   $\tilde \Sigma\in \mathcal V_n(M)$ whose support is a smooth embedded closed minimal hypersurface such that ${\bf F}(\tilde \Sigma, {\bf C}(\{\Phi_i\})_{i\in\N})<r, $
$${\bf L}(\Pi)=||\tilde \Sigma||(M)\text{ and index\,(support of $\tilde \Sigma$)}\leq k.$$
If we prove this  we obtain a sequence $\Sigma_j\in \mathcal V_n(M)$ whose  supports are smooth embedded closed minimal hypersurfaces of index less than or equal to $k$  and such that ${\bf F}(\Sigma_j, {\bf C}(\{\Phi_i\})_{i\in\N})<1/j$. By the monotonicity formula and the mass bound we get that the number of connected components of the support of $\Sigma_j$ and their multiplicities  are uniformly bounded. The Compactness Theorem proven by Sharp \cite{sharp} and the assumption that the metric is bumpy imply that there is a subsequence $\{m\}\subset \{j\}$ such that $\{\Sigma_m\}=\{\Sigma\}$.  This follows because otherwise we can find an embedded minimal surface of index $\leq k$ that is a nontrivial graphical limit away from finitely many points and with finite multiplicity of embedded minimal surfaces of index $\leq k$.   This must induce a nontrivial Jacobi field in the limit surface (see \cite{sharp}, Claim 6).
\smallskip

Denote by $\mathcal W$ the set of all  stationary integral varifolds  with mass $L$  whose support is a smooth embedded closed minimal hypersurface and by $\mathcal W(r)$ the set
$$\{V\in\mathcal W: {\bf F}(V, {\bf C}(\{\Phi_i\})_{i\in\N})\geq r\}.$$ 

\subsection{Lemma}{\em There exist $i_0\in \mathbb{N}$ and $\bar \varepsilon_0>0$ so that   $|\Phi_i|(X)\cap \overline{\bf B}^{\bf F}_{\bar \varepsilon_0}(\mathcal W(r))=\emptyset$ for all $i\geq i_0$.
}

\begin{proof}
If not there is a  subsequence $\{j\}\subset \{i\}$, and $x_j\in X$, $V_j\in\mathcal W(r)$ so that $$\lim_{j\to\infty}{\bf F}(|\Phi_j(x_j)|,V_j)=0.$$ Thus $\lim_{i\to\infty}||\Phi_j(x_j)||(M)=L$ and hence a subsequence of $|\Phi_j(x_j)|$ will converge to some $V \in {\bf C}(\{\Phi_i\})$.
 This is a contradiction. 
\end{proof}

Denote by $\mathcal W^{k+1}$ the collection of  elements in $\mathcal W$ whose support has  index greater than or equal to $(k+1)$.  Because the metric is bumpy this set is countable and thus we can write $\mathcal W^{k+1}\setminus  \overline{\bf B}^{\bf F}_{\bar \varepsilon_0}(\mathcal W(r))=\{\Sigma_1,\Sigma_2,\ldots\}$. 

Using Deformation Theorem A with $K=\overline{\bf B}^{\bf F}_{\bar \varepsilon_0}(\mathcal W(r))$ and $\Sigma=\Sigma_1$,  we find  $\bar \varepsilon_1>0$, $i_1\in\N$, and $\{\Phi^1_i\}_{i\in\N}$ in $\mathcal Z_n(M;{\bf F};G)$ so that
\begin{itemize}
\item $\Phi^1_i$ is homotopic to $\Phi_i$ in the ${\bf F}$-topology for all $i\in\N$;
\item $ {\bf L}(\{\Phi^1_i\}_{i\in\N})\leq L$;
\item for all $i\geq i_1$, $|\Phi^1_i|(X)\cap ({\overline{\bf B}}^{\bf F}_{\bar \varepsilon_1}(\Sigma_1)\cup  \overline{\bf B}^{\bf F}_{\bar \varepsilon_0}(\mathcal W(r)))=\emptyset$;
\item no $\Sigma_j$ belongs to $\partial {\overline{\bf B}}^{\bf F}_{\bar \varepsilon_1}(\Sigma_1)$.
\end{itemize}

We now look at $\Sigma_2$.  If $\Sigma_2 \notin {\overline{\bf B}}^{\bf F}_{\bar \varepsilon_1}(\Sigma_1)$, we apply Deformation Theorem A  and find  $\bar \varepsilon_2>0$, $i_2\in\N$, and $\{\Phi^2_i\}_{i\in\N}$ in $\mathcal Z_n(M;{\bf F};G)$ so that
\begin{itemize}
\item $\Phi^2_i$ is homotopic to $\Phi_i$ in the ${\bf F}$-topology for all $i\in\N$;
\item $ {\bf L}(\{\Phi^2_i\}_{i\in\N})\leq L$;
\item for all $i\geq i_2$, $|\Phi^2_i|(X)\cap ({\overline{\bf B}}^{\bf F}_{\bar \varepsilon_1}(\Sigma_1)\cup {\overline{\bf B}}^{\bf F}_{\bar \varepsilon_2}(\Sigma_2)\cup   \overline{\bf B}^{\bf F}_{\bar \varepsilon_0}(\mathcal W(r)))=\emptyset$;
\item no $\Sigma_j$ belongs to $\partial {\overline{\bf B}}^{\bf F}_{\bar \varepsilon_1}(\Sigma_1) \cup \partial {\overline{\bf B}}^{\bf F}_{\bar \varepsilon_2}(\Sigma_2)$.
\end{itemize}
If $\Sigma_2 \in {\bf B}^{\bf F}_{\bar \varepsilon_1}(\Sigma_1)$, we skip it and repeat the procedure  with $\Sigma_3$.

Proceeding inductively there are two possibilities. We can  find for all $l\in\N$ a sequence  $\{\Phi^l_i\}_{i\in\N}$, $\bar \varepsilon_l>0$, $i_l\in\N$, and $\Sigma_{j_l}\in  \mathcal W^{k+1}\setminus  \overline{\bf B}^{\bf F}_{\bar \varepsilon_0}(\mathcal W(r))$ for some $j_l\in\N$ so that
\begin{itemize}
\item $\Phi^l_i$ is homotopic to $\Phi_i$ in the ${\bf F}$-topology for all $i\in\N$;
\item $ {\bf L}(\{\Phi^l_i\}_{i\in\N})\leq L$;
\item for all $i\geq i_l$, $|\Phi^l_i|(X)\cap (\cup_{q=1}^l{\overline{\bf B}}^{\bf F}_{\bar \varepsilon_q}(\Sigma_{j_q})\cup   \overline{\bf B}^{\bf F}_{\bar \varepsilon_0}(\mathcal W(r)))=\emptyset$;
\item $\{\Sigma_1,\ldots,\Sigma_l\}\subset\cup_{q=1}^l{{\bf B}}^{\bf F}_{\bar \varepsilon_q}(\Sigma_{j_q})$;
\item no $\Sigma_j$ belongs to $\partial {\overline{\bf B}}^{\bf F}_{\bar \varepsilon_1}(\Sigma_1) \cup \cdots \cup \partial {\overline{\bf B}}^{\bf F}_{\bar \varepsilon_l}(\Sigma_{j_l})$.
\end{itemize}
Or the process ends in finitely many steps, which means we can find  some $m \in \mathbb{N}$, a sequence  $\{\Phi^m_i\}_{i\in\N}$, $\bar \varepsilon_1, \dots, \bar \varepsilon_m>0$, $i_m\in\N$, and $\Sigma_{j_1},\dots,\Sigma_{j_m}\in  \mathcal W^{k+1}\setminus  \overline{\bf B}^{\bf F}_{\bar \varepsilon_0}(\mathcal W(r))$ for some $j_l\in\N$ so that
\begin{itemize}
\item $\Phi^m_i$ is homotopic to $\Phi_i$ in the ${\bf F}$-topology for all $i\in\N$;
\item $ {\bf L}(\{\Phi^m_i\}_{i\in\N})\leq L$;
\item for all $i\geq i_m$, $|\Phi^l_i|(X)\cap (\cup_{q=1}^m{\overline{\bf B}}^{\bf F}_{\bar \varepsilon_q}(\Sigma_{j_q})\cup   \overline{\bf B}^{\bf F}_{\bar \varepsilon_0}(\mathcal W(r)))=\emptyset$;
\item $\{\Sigma_j: j\geq 1\}\subset\cup_{q=1}^m{{\bf B}}^{\bf F}_{\bar \varepsilon_q}(\Sigma_{j_q})$.
\end{itemize}

In the first case we choose an increasing sequence $p_l\geq i_l$ so that 
$$\sup_{x\in X}||\Phi^l_{p_l}||(M)\leq L+\frac 1 l,$$
and set $\Psi_l=\Phi^l_{p_l}$. In the second case we set $p_l=l$  and $\Psi_l=\Phi^m_l$.
The sequence $\{\Psi_l\}_{l\in\N}$ is such that 
\begin{itemize}
\item[(i)] $\Psi_l$ is homotopic to $\Phi_{p_l}$ in the ${\bf F}$-topology for all $l\in\N$;
\item[(ii)] $ {\bf L}(\{\Psi_l\}_{l\in\N})\leq L$;
\item[(ii)] ${\bf C}(\{\Psi_l\}_{l\in\N})\cap \left(\mathcal W^{k+1}\cup \mathcal W(r)\right)=\emptyset$.
\end{itemize}

The Min-max Theorem \ref{minmax.continuous.thm} implies $\mathcal W\setminus \left(\mathcal W^{k+1}\cup \mathcal W(r)\right)$ is nonempty and this is what we wanted to prove.


\end{proof}


\subsection{Remark}\label{xin.zhou.2} Again the proof gives more. One can choose $\Sigma$ to be $G$-almost minimizing as in Remark \ref{xin.zhou}.
\smallskip

We deduce Theorem \ref{theorem1} from the previous theorem.

\subsection*{Theorem \ref{theorem1}}{\em Let $(M^{n+1},g)$ be an $(n+1)$-dimensional closed Riemannian manifold,  with $3\leq (n+1)\leq 7$.  There exists an integral stationary varifold $\Sigma \in \mathcal V_n(M)$ with the following properties:
\begin{itemize}
\item $||\Sigma||(M)={\bf L}(\Pi),$
\item the support of $\Sigma$ is a smooth closed embedded hypersurface in $M$,
\item $index(support \, of \, \Sigma) \leq k$.
\end{itemize}
}

\begin{proof}
Consider a sequence of bumpy metrics $\{g_j\}_{j\in\N}$ tending smoothly  to   $g$ \cite{white3}. If $L_j$ is the width of $\Pi$ with respect to $g_j$,  the previous theorem gives    a  sequence $\{\Phi^j_i\}_{i\in\N}$ in $\Pi$ such that $L_j={\bf L}(\{\Phi_i^j\}_{i\in\N})$ and  $\Sigma_j\in {\bf C}(\{\Phi^j_i\}_{i\in\N})$ with support  a smooth embedded minimal hypersurface of index $\leq k$ such that $L_j=||\Sigma_j||(M)$ (with respect to the metric $g_j$). The sequence  $\{L_j\}_{j\in\N}$ tends to $L={\bf L}(\Pi)$ and so the result follows by using the Compactness Theorem A.6 of \cite{sharp} (for changing ambient metrics).
\end{proof}

Next we assume that $X$ is either $S^1$ or $[0,1]$, $G$ is either $\Z$ or $\Z_2$, 
$\Phi:X \rightarrow \mathcal Z_n(M^{n+1};{\bf F};G)$ is a continuous map, and  $\Pi$ is the class of all continuous maps $\Phi':X \rightarrow \mathcal Z_n(M^{n+1};{\bf F};G)$ such that $\Phi$ and $\Phi'$ are homotopic to each other in the flat topology.  If $X=[0,1]$, assume further that $\Phi'(0)=\Phi'(1)=0$ for all $\Phi'\in \Pi$ and that all homotopies are relative to $\partial [0,1]=\{0,1\}$.

 We first prove Theorem \ref{theorem2}  when the manifold $M$ does not contain one-sided embedded closed hypersurfaces. This is guaranteed if, for instance,  $H_n(M;\mathbb{Z}_2)=0$.

\subsection{Theorem}\label{index.1.theorem}{\em  We assume that $M$ contains no one-sided embedded closed hypersurfaces and  that $(M,g)$ is a bumpy metric. 
Let $\{\Phi_i\}_{i\in\N}$ be a  sequence in $\Pi$ such that $L={\bf L}(\{\Phi_i\}_{i\in\N})={\bf L}(\Pi).$ 
There is $\Sigma\in  {\bf C}(\{\Phi_i\}_{i\in\N})$ with support a smooth embedded closed minimal hypersurface
 such that $${\bf L}(\Pi)=||\Sigma||(M),\text{ index\,(support of $\Sigma$)}=1.$$
Moreover, the unstable component of $\Sigma$ has multiplicity one.}


\begin{proof}
Let $\mathcal W_j$ be the set of all  stationary integral varifolds  with mass $L$  whose support is a smooth embedded closed minimal hypersurface with index $j$, where $j=0,1$. Because the metric is bumpy, we obtain from \cite[Theorem 2.3]{sharp} (or more precisely its proof) that $\mathcal W_0\cup\mathcal W_1$ is finite. There can be at most one  unstable component in each  element in $\mathcal W_1$. We set  $\mathcal A$ to be those elements of $\mathcal W_0\cup\mathcal W_1$ that do not belong to ${\bf C}(\{\Phi_i\}_{i\in\N})$ and set $\mathcal U$ to be those elements of $\mathcal W_1\cap {\bf C}(\{\Phi_i\}_{i\in\N})$ for which its unstable component appears with multiplicity higher than one. We want to show that $\mathcal W_1\setminus (\mathcal A\cup \mathcal U)$ is non-empty.

There exist  $i_0\in\N$ and $\bar \varepsilon_0>0$ so that   $|\Phi_i|(X)\cap \overline{\bf B}^{\bf F}_{\bar \varepsilon_0}(\mathcal A)=\emptyset$ for all $i\geq i_0$.
Write $\mathcal U=\{\Sigma_1,\ldots, \Sigma_k\}$. Applying Deformation Theorem B to $\{\Phi_{i}\}_{i\in\N}$, with $K=\overline{\bf B}^{\bf F}_{\bar \varepsilon_0}(\mathcal A)$ and $\Sigma=\Sigma_1$, we find  $\bar \varepsilon_1>0$, $i_1\in\N$, and $\{\Phi^1_i\}_{i\in\N}$ in $\mathcal Z_n(M;{\bf F};G)$ so that
\begin{itemize}
\item $\Phi^1_i$ is homotopic to $\Phi_i$ in the flat topology for all $i\in\N$;
\item $ {\bf L}(\{\Phi^1_i\}_{i\in\N})\leq L$;
\item for all $i\geq i_1$, $|\Phi^1_i|(X)\cap {\overline{\bf B}}^{\bf F}_{\bar \varepsilon_1}(\{\Sigma_1\}\cup \mathcal A)=\emptyset$;
\item no $\Sigma_j$ belongs to $\partial  {\overline{\bf B}}^{\bf F}_{\bar \varepsilon_1}(\{\Sigma_1\}\cup \mathcal A)$.
\end{itemize}

 If $\Sigma_2\in {\bf C}(\{\Phi^1_i\}_{i\in\N})$, we apply Deformation Theorem B to a suitable subsequence of $\{\Phi^1_{i}\}_{i\in\N}$, with $K={\overline{\bf B}}^{\bf F}_{\bar \varepsilon_1}(\{\Sigma_1\}\cup \mathcal A)$ and $\Sigma=\Sigma_2$, to find  $\bar \varepsilon_2>0$, $i_2\in\N$, and $\{\Phi^2_i\}_{i\in\N}$ in $\mathcal Z_n(M;{\bf F};G)$ so that
\begin{itemize}
\item $\Phi^2_i$ is homotopic to $\Phi_i$ in the flat topology for all $i\in\N$;
\item $ {\bf L}(\{\Phi^2_i\}_{i\in\N})\leq L$;
\item for all $i\geq i_2$, $|\Phi^2_i|(X)\cap  {\overline{\bf B}}^{\bf F}_{\bar \varepsilon_2}(\{\Sigma_1,\Sigma_2\}\cup \mathcal A)=\emptyset$;
\item no $\Sigma_j$ belongs to $\partial  {\overline{\bf B}}^{\bf F}_{\bar \varepsilon_2}(\{\Sigma_1,\Sigma_2\}\cup \mathcal A)$.
\end{itemize}
If $\Sigma_2$ is not in ${\bf C}(\{\Phi^1_i\}_{i\in\N})$, then the above properties hold with $\Phi_i^2=\Phi_i^1$ and for some small $\bar \varepsilon_2$.

After a finite number of times we obtain $\{\Psi_{i}\}_{i\in\N}$ in $\mathcal Z_n(M;{\bf F};G)$, $i_k\in\N$, and $\delta>0$ so that 
\begin{itemize}
\item $\Psi_i$ is homotopic to $\Phi_i$ in the flat topology for all $i\in\N$;
\item $ {\bf L}(\{\Psi_i\}_{i\in\N})\leq L$;
\item for all $i\geq i_3$, $|\Psi_i|(X)\cap  {\overline{\bf B}}^{\bf F}_{\delta}(\mathcal U\cup \mathcal A)=\emptyset$.
\end{itemize}

Let $\hat \Psi_i$ be a pulled-tight sequence obtained from  $\{\Psi_i\}_{i\in\N}$. This means ${\bf C}(\{\hat \Psi_i\}) \subset {\bf C}(\{\Psi_i\})$  and every element of  ${\bf C}(\{\hat \Psi_i\})$ is stationary.  One can then find another $\delta'>0$ and $i_{k+1}\in\N$ so that  $|\hat \Psi_i|(X)\cap  {\overline{\bf B}}^{\bf F}_{\delta'}(\mathcal U\cup \mathcal A \cup \mathcal B)=\emptyset$ for all $i\geq i_{k+1}$, where $\mathcal B$ is the collection of elements of $\mathcal W_0$ that are not in ${\bf C}(\{\hat \Psi_i\}_{i\in\N})$.

Write $\mathcal W_0\setminus (\mathcal A \cup \mathcal B)=\{\tilde \Sigma_1,\ldots,\tilde \Sigma_l\}$. Every element in this set belongs to ${\bf C}(\{\hat \Psi_i\}_{i\in\N})$.
Then by Deformation Theorem C there exist $\xi>0$, $j_0\in\N$ so that  for all $i\geq j_0$ one can find   $\bar \Psi_i: X\rightarrow \mathcal Z_n(M;{\bf F};G)$ such that
\begin{itemize}
\item[(i)] $\bar \Psi_i$ is homotopic to $\hat \Psi_i$  in the flat topology;
\item[(ii)] ${\bf L}(\bar \Psi_i\})\leq L$;
\item[(iii)]  $${\bf \Lambda}(\{\bar \Psi_i\})\subset \left({\bf \Lambda}(\{\hat \Psi_i\})\setminus \cup_{q=1}^l {\bf B}_\xi^{\bf F}(\tilde \Sigma_q)\right) \cup \left(\mathcal V_n(M) \setminus {\bf B}_\xi^{\bf F}(\Gamma)\right),$$
where $\Gamma$ is the collection of all stationary integral varifolds  $V\in \mathcal V_n(M)$ with $L=||V||(M)$.
\end{itemize}

In particular, ${\bf C}(\{\bar \Psi_i\}_{i\in\N})\cap (\mathcal U\cup  \mathcal W_0\cup \mathcal A)=\emptyset$. Theorem \ref{index.k.theorem} applied to $\{\bar \Psi_i\}$ implies that $(\mathcal W_0\cup\mathcal W_1)\setminus(\mathcal U\cup  \mathcal W_0\cup \mathcal A)$ is nonempty and this is what we wanted to prove.


\end{proof}

The next proposition is similar to arguments of X. Zhou \cite{zhou}, who considered $G=\mathbb{Z}$ and the multiplicity of non-orientable components. 

\subsection{Proposition}\label{even.multiplicity}{\em Suppose $\Sigma\in \mathcal V_n(M)$ is stationary and $G$-almost minimizing of boundary type (see Remark \ref{xin.zhou}) in annuli. Write $$\Sigma=m_1|S_1|+\cdots + m_{P}|S_{P}|,$$ where $m_i\in \mathbb{N}$ and $S_1,\dots, S_P$ are the connected components of ${\rm supp\,} \Sigma$ (smooth embedded closed minimal hypersurface). For any one-sided component $S_i$, $m_i$ is an even number.}

\begin{proof}
As in the proof of Proposition 6.1 of \cite{zhou}, we can find $p \in S_i$, $s>0$, a sequence $T_j\in \mathcal Z_n(M;G)$, $T_j=\partial U_j$, such that ${\bf F}(|T_j|,\Sigma) \rightarrow 0$ as $j\rightarrow \infty$ and ${\rm supp\,} T_j\llcorner B_s(p)$ is a smooth embedded stable minimal hypersurface.  By compactness, we can assume $T_j\rightarrow T=\partial U$ in the flat topology. For some tubular neighborhood of $S_i$, 
$T \llcorner \Omega \in \mathcal Z_n(M, M \setminus \Omega;G)$ is a relative boundary. But ${\rm supp\,}(T \llcorner \Omega) \subset S_i$. Since $S_i$ is one-sided, the Constancy Theorem implies $T \llcorner \Omega=k_i S_i$ with $k_i$ even (if $G=\mathbb{Z}$ and $S_i$ is non-orientable, necessarily $k_i=0$). By stability and Schoen-Simon regularity theory \cite{schoen-simon}, the convergence of $T_j$ to $\Sigma$ is smooth and locally graphical  near $p$. Since $|a|\equiv a {\rm \, mod\,}2$ for $a\in G$, we have $m_i \equiv k_i {\rm \, mod\,}2$.
\end{proof}

We can now prove Theorem \ref{theorem2} and Corollary \ref{corollary1}.

\subsection*{Theorem  \ref{theorem2}}\label{index.onesided}{\em  Suppose $(M,g)$ is a bumpy metric, and $X=[0,1]$. We assume  $\{\Phi_i\}_{i\in\N}$ is a  sequence in $\Pi$ such that $L={\bf L}(\{\Phi_i\}_{i\in\N})={\bf L}(\Pi).$ 
There exists $\Sigma\in  {\bf C}(\{\Phi_i\}_{i\in\N})$ with support a smooth embedded closed minimal hypersurface 
 such that ${\bf L}(\Pi)=||\Sigma||(M)$ and either:
 \begin{itemize}
 \item[(a)] $\text{ index\,(support of $\Sigma$)}=1,$ and the unstable component of $\Sigma$ has multiplicity one;
 \item[(b)] $\text{ index\,(support of $\Sigma$)}=0,$ and at least one component of $\Sigma$ is one-sided, has even multiplicity and  unstable double cover.
 \end{itemize}
}

\begin{proof}
 Let $\bar \Sigma$ be an embedded one-sided closed minimal hypersurface. Since the metric is bumpy, the double cover is either unstable or
strictly stable. If it is strictly stable, one can construct $P_t:\Omega_{r_0} \rightarrow \Omega_{r_0}$ just like in Proposition  \ref{projection}, where
$\Omega_{r_0}$ is a neighborhood of $\bar \Sigma$ in $M$. This means that embedded one-sided closed minimal hypersurfaces with strictly stable double covers can be treated just like strictly stable two-sided minimal hypersurfaces if one follows the proof of Theorem \ref{index.1.theorem} and Deformation Lemma C. With the appropriate modifications, we get the existence of $\Sigma\in  {\bf C}(\{\Phi_i\}_{i\in\N})$ with support a smooth embedded closed minimal hypersurface 
 such that ${\bf L}(\Pi)=||\Sigma||(M)$ and either:
 \begin{itemize}
 \item[(a)] $\text{ index\,(support of $\Sigma$)}=1,$ and the unstable component of $\Sigma$ has multiplicity one;
 \item[(b)] $\text{ index\,(support of $\Sigma$)}=0,$ and at least one component of $\Sigma$ is one-sided, and  has unstable double cover.
 \end{itemize}
Since  $\Sigma$ can be chosen to be $G$-almost minimizing of boundary type in annuli, any one-sided component of it has even multiplicity by Proposition \ref{even.multiplicity}. This proves Theorem \ref{index.onesided}.

\end{proof}

\subsection*{Corollary \ref{corollary1}}{\em  Suppose $(M,g)$ contains no one-sided embedded hypersurface, and $X=[0,1]$.
There exists an integral stationary varifold $\Sigma \in \mathcal V_n(M)$ with $||\Sigma||(M)={\bf L}(\Pi)$ and satisfying  the following properties:
\begin{itemize}
\item[(a)] the support of $\Sigma$ is a smooth closed embedded minimal hypersurface in $M$;
\item[(b)]  we have \begin{eqnarray*}
index(support \, of \, \Sigma)  &\leq& 1\\
&\leq& index(support \, of \, \Sigma) +nullity(support \, of \, \Sigma);
\end{eqnarray*}
\item[(c)] any unstable component (necessarily unique)  of $\Sigma$ must have  multiplicity one;
\item[(d)] if the metric $g$ is bumpy and $\{\Phi_i\}_{i\in\N}$ is a min-max sequence, $\Sigma$ can be chosen so that for some subsequence 
$i_j \rightarrow \infty$  there is some $x_{i_j}\in X$ such that $\Phi_{i_j}(x_{i_j})$ converges in varifold sense to $\Sigma$.
\end{itemize}
}

\begin{proof}
For bumpy metrics the result follows at once from Theorem \ref{theorem2}. For a general metric $g$,  consider a sequence of bumpy metrics $\{g_j\}_{j\in\N}$ tending smoothly  to   $g$ (by \cite{white3}). If $L_j$ is the width of $\Pi$ with respect to $g_j$,  we obtain  a sequence   $\{\Sigma_j\}_{j\in\N}\in \mathcal V_n(M)$ satisfying
\begin{itemize}
\item the support of $\Sigma_j$ is a smooth closed embedded minimal hypersurface in $M$;
\item$index(support \, of \, \Sigma_j)=1;$
\item the  unstable component   of $\Sigma_j$ has  multiplicity one;
\item $L_j=||\Sigma_j||(M)$ (with respect to the metric $g_j$) for all $j\in\N$.
\end{itemize}

The sequence  $\{L_j\}_{j\in\N}$ tends to $L={\bf L}(\Pi)$ and the Compactness Theorem A.6 of Sharp \cite{sharp} (for changing ambient metrics) implies that, after passing to a subsequence, $\{\Sigma_j\}_{j\in\N}$ converges in the varifold sense to an integral stationary varifold  $\Sigma\in\mathcal V_n(M)$ satisfying (a) and such that $||\Sigma||(M)={\bf L}(\Pi)$. Moreover, it also follows from the same result  that (after passing to a further subsequence):
\begin{itemize}
\item the union of the stable components of $\Sigma_j$  with their respective multiplicities converges in the varifold sense to  an integral stationary varifold $\Sigma_s$  with support a smooth stable closed embedded minimal hypersurface in $M$;
\item  the unstable component of $\Sigma_j$ with multiplicity one  converges in the varifold sense  to an integral stationary varifold $\Sigma_u$ with support a smooth closed embedded minimal hypersurface in $M$.
\end{itemize}
  Either the convergence to $\Sigma_u$ is smooth and graphical, in which case $\Sigma_u$ has multiplicity one and satisfies (b) and (c), or the convergence is not smooth. In the latter case  $\Sigma_u$ must be  stable with  nullity (see \cite{sharp}, Claim 6), and therefore  also satisfies (b).
Since $\Sigma=\Sigma_s+\Sigma_u$, we are done. 
\end{proof}

\appendix

\section{Various Interpolation Results}\label{A.section}


\subsection{Proposition}\label{mass.close}\textit{There exist constants $\alpha=\alpha(M)>0$ and $A_1=A_1(M)>0$ such that, if $\Phi_1,\Phi_2:[0,1] \rightarrow  \mathcal{Z}_n(M;\mathcal F;G)$ are continuous maps in the flat topology with $\Phi_1(0)=\Phi_2(0)$, $\Phi_1(1)=\Phi_2(1)$ and 
$$
\sup\{{\bf M}(\Phi_1(x)-\Phi_2(x)): x\in [0,1]\} < \alpha,
$$
then there exists a homotopy $H$ of $\Phi_1$ and $\Phi_2$, continuous in the flat topology and relative to $\partial [0,1]=\{0,1\}$, such that
$$
\sup\{{\bf M}(H(x,t)-\Phi_i(x)): x\in [0,1], t\in [0,1]\} < A_1\alpha,
$$
where $i=1,2$.}

\begin{proof}
Let $\Psi(x)=\Phi_2(x)-\Phi_1(x)$. Note that $\mathcal F(\Psi(x))\leq {\bf M}(\Psi(x))$, $\Psi(0)=\Psi(1)=0$. By Theorem 8.2 of Almgren \cite{almgren}, if $\alpha$ is sufficiently small there exists a homotopy $h(x,t)$ between $\Psi$ and $0$ such that 
$$
\sup\{{\bf M}(h(x,t)): x \in [0,1], t\in [0,1]\}\leq \mu_1 \sup\{{\bf M}(\Psi(x)): x\in [0,1]\}.
$$
The result follows by defining $H(x,t)= \Phi_2(x)-h(x,t)$.
\end{proof}

A related result to the next proposition was also proven by Zhou \cite{zhou2}.

\subsection{Proposition}\label{f.close} \textit{Given $\delta>0$, $L>0$, let $\mathcal S \subset \mathcal Z_n(M;{\bf F}; G)$ be a compact set
with ${\bf M}(R)\leq 2L$ for every $R\in \mathcal S$. Then there exists $\eta=\eta(\mathcal S,\delta, L)<\delta$ so that for all  $S\in \mathcal Z_n(M;G)$ with  
$${\bf M}(S)\leq 2L\quad\mbox{and}\quad\mathcal F(S,R)<\eta\quad\mbox{for some }R\in \mathcal S,$$
 there exists  $H:[0,1]\rightarrow  \mathcal{Z}_n(M;{\bf M};G)$ continuous in the mass topology and such that
\begin{itemize}
\item $H(0)=S$, $H(1)=R$, and $H([0,1])\subset {\bf B}^{\mathcal F}_{\delta}(R); $
\item ${\bf M}(H(x))\leq {\bf M}(S)+\delta$ for all $0\leq x\leq 1$.
\end{itemize}
Moreover, if ${\rm supp}\, R \subset W \subset \subset U$ for all $R\in \mathcal{S}$ and some open sets $U,W$, we can choose $\eta=\eta(\mathcal S,\delta, L,W, U)<\delta$ such that if $S$ is as above and ${\rm supp}\, S\subset W$ we can choose $H$ with ${\rm supp}\, H(t)\subset U$ for every $t\in [0,1]$.}

\begin{proof}

The next two auxiliary lemmas were essentially proven in \cite{pitts}.

\subsection{Lemma}\label{lemma3.7}{\em Let  $W\subset \subset U$ be  open subsets of $M$.  If $V\in \mathcal V_n(M), T, T_1,\ldots,$ are elements in $\mathcal Z_n(W;G)$ with 
$$\lim_{i\to\infty}(\mathcal F(T,T_i)+ {\bf F}(V,|T_i|))=0,$$
then for every $\bar \delta>0$ there exist  a sequence $S_1,S_2,\ldots,$ in $\mathcal Z_n(U;G)$ such that  
$\lim_{i\to\infty} {\bf F}(T,S_i)=0$
and the following property holds:

For each $i\in\N$ there is $m_i\in\N$ and a finite sequence $R^i_0,\ldots,R^i_{m_i}$ in $\mathcal Z_n(U;G)$ such that $R^i_0=T_i$, $R^i_{m_i}=S_i$,
$${\bf M}(R^i_j)\leq {\bf M}(T_i)+\bar\delta,\quad {\bf M}(R^i_{j}-R^i_{j-1})\leq \bar\delta \quad\mbox{for all }j=1,\ldots,m_i$$
and $\lim_{i\to\infty}\sup_{0\leq j_i\leq m_i}\mathcal F(R^i_{j_i},T)=0.$}
\begin{proof}
In Lemma 3.7 of \cite{pitts} (using $T^*=0$ and $K=\bar V$ in his notation), all properties are proven except for the fact that $R^i_{j_i}$ tends to $T$ uniformly in the flat topology as $i\to\infty$. While Pitts did not state it, this property follows from the proof of Lemma 3.7, as we now explain. We refer to the notation used by Pitts. The currents  $R^i_j$ are constructed either in  \cite[page 117]{pitts}, and one sees at once that $\mathcal F(R^i_j,T)\leq {\bf M}(Q_i)=\mathcal F(T_i,T)$ for all $i\in\N$ and $0\leq j\leq m_i$,   or in  \cite[page 121]{pitts} and one sees that, for all $0\leq j\leq m_i$, $R^i_j$ coincides with $T_i$ on $M\setminus B_{r_i}(p_i)$, where $r_i$ tends to zero and $p_i$ tends to $p$, as $i\to\infty$.  In both cases  $R^i_{j_i}$ tends to $T$ uniformly in $j_i$as $i\to\infty$.

\end{proof}

\subsection{Lemma}\label{lemma3.8}{\em Let  $W\subset \subset U$ be  open subsets of $M$.  Given $\bar L>0$ and $\bar \delta>0$, there exists $\varepsilon=\varepsilon(W,U,\bar L,\bar \delta)>0$ such that whenever $S_1,S_2\in \mathcal Z_n(W;G)$ have $\mathcal F(S_1,S_2)\leq \varepsilon$ and ${\bf M}(S_1)\leq \bar L$, ${\bf M}(S_2)\leq \bar L$, there is a finite sequence $T_0,\ldots,T_m$ in $\mathcal Z_n(U;G)$ such that $T_0=S_1$, $T_{m}=S_2$,
$${\bf M}(T_j)\leq \bar L+\bar\delta,\quad {\bf M}(T_{j}-T_{j-1})\leq \bar\delta \quad\mbox{for all }j=1,\ldots,m$$
and $\mathcal F(T_j,S_1)<\bar \delta$ for all $0\leq j\leq m$.}
\begin{proof}
In Lemma 3.8 of \cite{pitts} (using $T=0$ and $K=\bar W$ in his notation), all properties are proven except for the fact that $\mathcal F(T_j,S_1)<\bar \delta$ for all $0\leq j\leq m$. To prove this last property it suffices to do the following modification of the proof of Lemma 3.8 of \cite{pitts} (we now use his notation): On page 122 of \cite{pitts}, use Lemma \ref{lemma3.7} instead of Lemma 3.7. The whole argument carries over but one needs to check that  the currents $T_j$ defined at the bottom of page 123 can be made arbitrarily close to $S_i$ in the flat metric and this is indeed the case because $\mathcal F(S_0,T_j)\leq {\bf M}(Q)=\mathcal F(S_0,S_i)$.
\end{proof}

Let $C_0=C_0(M,1)>0$ and $\delta_0=\delta_0(M)>0$ be as in Theorem 3.10 of \cite{marques-neves-infinitely}. We can assume $\delta < \delta_0$ and set $\bar \delta=\delta/(3+C_0).$

We argue by contradiction and suppose there is a sequence $S_i\in \mathcal Z_n(M;G)$, ${\bf M}(S_i)\leq 2L$, such that $\mathcal F(S_i,R_i)<1/i$ for some $R_i\in \mathcal S$ and for which the conclusion of the proposition does not hold. After passing to a subsequence if necessary, we find $V\in\mathcal V_n(M)$ and $R\in \mathcal S$ such that
$$\lim_{i\to\infty}({\bf F}(R_i,R)+ {\bf F}(|S_i|,V))=0\quad\mbox{and}\quad {\bf M}(R)\leq ||V||(M).$$
Set $\bar L=||V||(M)+\bar\delta$ and consider $\varepsilon=\varepsilon(M,M,\bar L,\bar\delta)$ given by Lemma \ref{lemma3.8}. Choose $i$ sufficiently large so that 
$$||S_i||(M)\leq \bar L\leq ||S_i||(M)+2\bar\delta, \,\, {\bf M}(R_i)\leq\bar L,\mbox{ and } \mathcal F(S_i,R_i)<\varepsilon.$$
Applying Lemma \ref{lemma3.8} to $S_i, R_i$ and using Theorem 3.10 of \cite{marques-neves-infinitely} we obtain $H:[0,1] \rightarrow  \mathcal Z_n(M;{\bf M}; G)$ continuous in the mass topology such that $H(0)=S_i$, $H(1)=R_i$,
 $$\sup_{x\in [0,1]}{\bf M}(H(x)) \leq  \bar L +(C_0+1)\bar\delta\leq {\bf M}(S_i) + (C_0+3) \bar\delta\leq {\bf M}(S_i)+\delta,$$ and 
 $H([0,1])\subset {\bf B}^{\mathcal F}_{(C_0+1)\bar \delta}(R_i)\subset {\bf B}^{\mathcal F}_{\delta}(R_i)$. This is a contradiction and finishes the first part of the proof of the Proposition. 
 
 To prove the last statement  one needs to localize Theorem 3.10 of \cite{marques-neves-infinitely} so that the deformations take place near $U$. This follows because Almgren interpolation results and the proof of Theorem 3.10 of \cite{marques-neves-infinitely} use a triangulation of $M$ and so one chooses a triangulation sufficiently fine with respect to $d(U,M\setminus W)$ (see proof of Theorem 14.1 of \cite{marques-neves}).

\end{proof}

\subsection{Proposition}\label{F.approximation} {\em Let $Y$ be a cubical subcomplex of some $I(m,l)$ and let $\Phi:Y\rightarrow \mathcal{Z}_n(M;{\bf F};G)$ be
a continuous map in the ${\bf F}$-topology. 
There exist a sequence of maps continuous in the mass topology
$$\Phi_i:Y \rightarrow \mathcal{Z}_n(M;{\bf M};G),\quad  j\in\N, $$
homotopic to $\Phi$ in the flat topology, and such that 
 $$\sup\{ {\bf F}(\Phi_i(y),\Phi(y)): y\in Y\}\rightarrow 0$$
 as $i\rightarrow \infty$.
 If $Z$ is a subcomplex of $Y$  such that $\Phi_{|Z}\equiv 0$, then the sequence $\{\Phi_i\}_{i\in\N}$ can be chosen 
so that $(\Phi_{i})_{|Z}\equiv 0$ for all $i\in\N$ and the homotopies can be taken relative to $Z$.
}

\begin{proof}
Given $\Phi: Y \rightarrow \mathcal Z_n(M;{\bf F};G)$, Theorem 3.9 of \cite{marques-neves-infinitely} (the analogous statement holds for $G=\mathbb{Z}$) implies there exist a sequence of maps
$$\phi_i:Y(k_i)_0 \rightarrow \mathcal{Z}_n(M;G),$$
with $k_i<k_{i+1}$, and a
sequence of positive numbers $\{\delta_i\}_{i\in\N}$ converging to zero such that
\begin{itemize}
\item[(i)] $$S=\{\phi_i\}_{i\in\N}$$ is an $(Y,{\bf M})$-homotopy sequence of mappings into $\mathcal{Z}_n(M;{\bf M};G)$ with ${\bf f}(\phi_i)<\delta_i$;
\item[(ii)] $$\sup\{\mathcal F(\phi_i(x)-\Phi(x)): x\in Y(k_i)_0\}\leq \delta_i;$$
\item[(iii)]$$\sup\{{\bf M}(\phi_i(x)): x\in X(k_i)_0\}\leq \sup\{{\bf M}(\Phi(x)): x\in Y\}+\delta_i.$$
\end{itemize}
We can also guarantee that (see Theorem 13.1 (i) of \cite{marques-neves} for the case $Y=I^m$)
\begin{itemize}
\item[(iv)] for some $l_i\rightarrow \infty$  and every $y\in X(k_i)_0$,
$$
{\bf M}(\phi_i(y))\leq \sup \{{\bf M}(\Phi(x)):\alpha \in X(l_i), x,y\in \alpha\} + \delta_i.
$$
\end{itemize}

Because $x\in Y \mapsto {\bf M}(\Phi(x))$ is continuous, we get from property (iv) above that for every $y\in X(k_i)_0$,
$$
{\bf M}(\phi_i(y)) \leq {\bf M}(\Phi(y))+\eta_i
$$
with $\eta_i\rightarrow 0$ as $i \rightarrow \infty$.
If we apply Lemma 4.1 of \cite{marques-neves} with $\mathcal S=\Phi(X)$,  property [(ii)] above implies
\begin{equation}\label{appendix.convergence}
\sup\{{\bf F}(\phi_i(x),\Phi(x)): x\in X(k_i)_0\} \rightarrow 0
\end{equation}
as $i \rightarrow \infty$.  

For sufficiently large $i$, by Theorem 3.10 of \cite{marques-neves-infinitely} we can take the Almgren extension of $\phi_i$, $\Phi_i: Y \rightarrow \mathcal Z_n(M;{\bf M};G)$. By Corollary 3.12 of 
 \cite{marques-neves-infinitely}, $\Phi_i$ and $\Phi$ are homotopic to each other in the flat topology. By Theorem 3.10 (iii) \cite{marques-neves-infinitely}, (\ref{appendix.convergence}) and since $\Phi$ is continuous in the $F$-metric, we have that
  $$\sup\{ {\bf F}(\Phi_i(y),\Phi(y)): y\in Y\}\rightarrow 0$$
  as $i\rightarrow \infty$.
  
If  $\Phi_{|Z}\equiv 0$, then naturally $(\Phi_i)_{|Z}\equiv 0$ and the fact that the homotopies can be taken relative to $Z$ follows from the
proof of Proposition 3.5 in the Appendix of \cite{marques-neves-infinitely}.
\end{proof}

\bibliographystyle{amsbook}

\begin{thebibliography}{99}


\bibitem{agol-marques-neves}
Agol, I., Marques, F. C., Neves, A., \textit{Min-max theory and the energy of links,} \emph{arXiv:1205.0825 [math.GT]}(2012), 
to appear in Journal of the Amer. Math. Soc.  	

\bibitem{nicolau} N. Aiex, \textit{The width of Ellipsoids,} preprint.

\bibitem{allard}
W. K. Allard, 
\textit{On the first variation of a varifold,}
Ann. of Math. (2) \textbf{95} (1972), 417--491. 




\bibitem{almgren} 
F. Almgren, \textit{The homotopy groups of the integral cycle groups,} Topology  (1962), 257--299. 

\bibitem{almgren-varifolds}
F. Almgren, \textit{The theory of varifolds,} Mimeographed notes, Princeton (1965).




\bibitem{birkhoff} G. Birkhoff, \textit{Dynamical systems with two degrees of freedom,} Trans. Amer. Math. Soc. 18 (1917), 199--300.

\bibitem{panov} V. Buchstaber and T. Panov, 
\textit{Torus actions and their applications in topology and combinatorics,} 
University Lecture Series, 24. American Mathematical Society, Providence, RI, 2002. viii+144 pp.

\bibitem{buzano-sharp}
R. Buzano and B. Sharp, \textit{Qualitative and quantitative estimates for minimal hypersurfaces with bounded index and area,}
 arXiv:1512.01047v2 [math.DG]


\bibitem{carlotto}
A. Carlotto, \textit{Generic finiteness of minimal surfaces with bounded Morse index,} arXiv:1509.07101v3 [math.DG]


\bibitem{chodosh-ketover-maximo}
O. Chodosh, D. Ketover and D. Maximo, \textit{Minimal hypersurfaces with bounded index,}  arXiv:1509.06724v3 [math.DG]

\bibitem{colding-delellis} 
T. Colding and C. De Lellis, \textit{The min-max construction of minimal surfaces,}   Surveys in  Differential Geometry VIII , International Press,  (2003),   75--107.




\bibitem{delellis-genus} 
C. De Lellis and F. Pellandini, \textit{Genus bounds for minimal surfaces arising from min-max constructions,} 
J. Reine Angew. Math. 644 (2010), 47--99.


\bibitem{federer}
H. Federer, \textit{Geometric measure theory,} 
Die Grundlehren der mathematischen Wissenschaften, Band 153 Springer-Verlag New York Inc., New York 1969. 














\bibitem{guaraco}
M. Guaraco, \textit{Min-max for phase transitions and the existence of embedded minimal hypersurfaces}, arXiv:1505.06698v2 [math.DG]







\bibitem{ketover}
D. Ketover, \textit{Degeneration of Min-Max Sequences in 3-manifolds,} 	arXiv:1312.2666 [math.DG] (2013)

\bibitem{ketover-marques-neves}
D. Ketover, F. C. Marques and A. Neves, \textit{The catenoid estimate and its geometric applications,} preprint

\bibitem{ketover-zhou}
D. Ketover and X. Zhou, \textit{Entropy of closed surfaces and min-max theory,}
	arXiv:1509.06238 [math.DG] (2015)







\bibitem{li-zhou}
H. Li and X. Zhou, \textit{Existence of minimal surfaces of arbitrary large Morse index}, 	arXiv:1504.00970 [math.DG]




\bibitem{marques-neves-duke}
Marques, F. C.,  Neves, A.,
\textit{Rigidity of min-max minimal spheres in three-manifolds,}
Duke Math. J. \textbf{161} (2012), no. 14, 2725--2752. 

\bibitem{marques-neves} Marques, F. C.,  Neves A., \textit{Min-max theory and the Willmore conjecture,} Ann. of Math. \textbf{179} 2 (2014), 683--782.

\bibitem{marques-neves-infinitely}
Marques, F. C., Neves, A.,
\textit{Existence of infinitely many minimal hypersurfaces in positive Ricci curvature,}
arXiv:1311.6501 [math.DG] (2013).

\bibitem{marques-neves-cdm}
Marques, F. C., Neves, A.,
\textit{Applications of Almgren-Pitts min-max theory,} Current developments in mathematics 2013, 1--71, Int. Press, Somerville, MA, 2014. 

\bibitem{mazet-rosenberg}
L. Mazet and H. Rosenberg, \textit{Minimal hypersurfaces of least area},  arXiv:1503.02938v2 [math.DG]


\bibitem{montezuma} R. Montezuma, \textit{Min-max minimal hypersurfaces in non-compact manifolds}, 	arXiv:1405.3712 (2014), to appear in J. Differ. Geom.










\bibitem{pitts} 
J. Pitts, \textit{Existence and regularity of minimal surfaces on Riemannian manifolds,} Mathematical Notes 27, Princeton University Press, Princeton, (1981).

\bibitem{pitts-rubinstein}
J. Pitts, J. H. Rubinstein, J. H.,
\textit{Existence of minimal surfaces of bounded topological type in three-manifolds,} Proc. Centre Math. Anal. Austral. Nat. Univ., 10, Austral. Nat. Univ., Canberra, 1986.














\bibitem{schoen-simon} 
R. Schoen and L. Simon, \textit{Regularity of stable minimal hypersurfaces,} 
Comm. Pure Appl. Math. 34 (1981), 741--797. 
\bibitem{sharp}
B. Sharp, {\em Compactness of minimal hypersurfaces with bounded index,} preprint.


\bibitem{smith} 
Smith, On the existence of embedded minimal 2Ð-spheres in the 3--sphere, endowed with an arbitrary Riemannian metric,  \emph{Ph.D. thesis, supervisor L. Simon, University of Melbourne} (1982).





\bibitem{song}
A. Song, \textit{Embeddedness of least area minimal hypersurfaces,} 	arXiv:1511.02844 [math.DG]


 



\bibitem{urbano}
Urbano, F.,
\textit{Minimal surfaces with low index in the three-dimensional sphere,}
Proc. Amer. Math. Soc. \textbf{108} (1990), 989--992.

 
 

\bibitem{white2}
B. White, \textit{The space of minimal submanifolds for varying Riemannian metrics,} Indiana Univ. Math. J. 40 (1991), 161--200.

\bibitem{white3}
B. White, \textit{On the Bumpy Metrics Theorem for Minimal Submanifolds,} 	arXiv:1503.01803 [math.DG] (2015)







\bibitem{yau1} S.-T. Yau \textit{
Problem section.} Seminar on Differential Geometry, pp. 669--706, 
Ann. of Math. Stud., 102, Princeton Univ. Press, Princeton, N.J., 1982.

\bibitem{zhou}
X. Zhou,
\textit{Min-max minimal hypersurface in $(M^{n+1},g)$ with $Ric>0$ and $2\leq n\leq 6$.}
J. Differential Geom. 100 (2015), no. 1, 129--160. 

\bibitem{zhou2}
X. Zhou,
\textit{Min-max hypersurface in manifold of positive Ricci curvature,}
	arXiv:1504.00966 [math.DG] (2015)

\end{thebibliography}

\end{document}